\documentclass[aop,MSNbibl,seceqn,dvips]{arximspdf}
\usepackage{eufrak}

\doi{10.1214/12-AOP826} 
\volume{42}
\issue{2}
\pubyear{2014}
\firstpage{497}
\lastpage{526}

\makeatletter

\newcommand{\rrVert}{\Vert}
\newcommand{\rrvert}{\vert}
\newcommand{\llVert}{\Vert}
\newcommand{\llvert}{\vert}
\newcommand{\eqref}[1]{(\ref{#1})}
\newtheorem{prop}{Proposition}[section]
\newtheorem{cor}[prop]{Corollary}
\newtheorem{lemme}[prop]{Lemma}
\newtheorem{lemma}[prop]{Lemma}
\newproclaim{rem}[prop]{Remark}
\newproclaim{remark}[prop]{Remark}
\newproclaim{rems}[prop]{Remarks}
\newtheorem{thmm}[prop]{Theorem}
\newproclaim{defi}[prop]{Definition}
\newproclaim{example}[prop]{Example}
\newproclaim{examples}[prop]{Examples}

\newcommand{\N}{\mathbb{N}}
\newcommand{\Z}{\mathbb{Z}}
\newcommand{\R}{\mathbb{R}}

\newcommand{\norm}[2][]{#1\llVert  #2#1\rrVert }

\def\HH{\EuFrak H}

\def\1{{\mathbf{1}}}

\def\1{{\mathbf{1}}}
\makeatother

\begin{document}
\begin{frontmatter}

\title{Asymptotic independence of multiple Wiener--It\^o integrals and
the resulting limit laws}
\runtitle{Asymptotic independence and limit laws}

\begin{aug}
\author[A]{\fnms{Ivan} \snm{Nourdin}\corref{}\ead[label=e1]{inourdin@gmail.com}\ead[label=u1,url]{http://www.iecn.u-nancy.fr/\textasciitilde nourdin/}\thanksref{t1}}
\and
\author[B]{\fnms{Jan} \snm{Rosi\'nski}\ead[label=e2]{rosinski@math.utk.edu}\ead[label=u2,url]{http://www.math.utk.edu/\textasciitilde rosinski/}}
\runauthor{I. Nourdin and J. Rosi\'nski}
\affiliation{Universit\'{e} de Lorraine and University of Tennessee}
\thankstext{t1}{Supported in part by ANR Grants ANR-09-BLAN-0114
and ANR-10-BLAN-0121 at Universit\'e de Lorraine.}
\address[A]{Institut de Math\'{e}matiques \'Elie Cartan\\
Universit\'{e} de Lorraine\\
B.P. 70239\\
54506 Vandoeuvre-l\`es-Nancy Cedex\\
France\\
\printead{e1}\\
\printead{u1}} 
\address[B]{Department of Mathematics\\
University of Tennessee\\
Knoxville, Tennessee 37996\\
USA\\
\printead{e2}\\
\printead{u2}}
\end{aug}

\received{\smonth{1} \syear{2012}}
\revised{\smonth{10} \syear{2012}}

%
\begin{abstract}
We characterize the asymptotic independence between blocks consisting
of multiple Wiener--It\^o integrals.
As a consequence of this characterization, we derive the celebrated
fourth moment theorem of Nualart and Peccati, its multidimensional
extension and other related results on the multivariate convergence of
multiple Wiener--It\^o integrals, that involve Gaussian and non
Gaussian limits. We give applications to the study of the asymptotic
behavior of functions of short and long-range dependent stationary
Gaussian time series and establish the asymptotic independence for
discrete non-Gaussian chaoses.
\end{abstract}

%
\pdfkeywords{60F05, 60G15, 60H05, 60H07, Multiple Wiener-It\^o integral, multiplication formula, limit theorems}
\end{frontmatter}

\section{Introduction}
Let $B=(B_t)_{t\in\R_+}$ be a standard one-dimensional Brownian motion,
$q\geq1$ be an integer and let $f$ be a symmetric element of $L^2(\R
^q_+)$. Denote by $I_q(f)$ the $q$-tuple Wiener--It\^o integral of $f$
with respect to $B$. It is well known that multiple Wiener--It\^o
integrals of different orders are uncorrelated but not necessarily
independent. In an important paper \cite{UZ}, \"Ust\"unel and Zakai
gave the following characterization of the independence of multiple
Wiener--It\^o integrals.

\begin{thmm}[(\"Ust\"unel--Zakai)]\label{UZthmIntro}
Let $p,q\geq1$ be integers and let $f\in L^2(\R^p_+)$ and $g\in L^2(\R
^q_+)$ be symmetric. Then, random variables $I_p(f)$ and $I_q(g)$ are
independent if and only if
%
\begin{eqnarray}
\label{eqUZ} &&\int_{\R^{p+q-2}_+} \biggl\llvert \int
_{\R_+} f(x_1,\ldots,x_{p-1},u)
\nonumber
\\[-8pt]
\\[-8pt]
\nonumber
&&\hspace*{54pt}{}\times g(x_{p+1},
\ldots,x_{p+q-2},u) \,du \biggr\rrvert ^2 \,dx_1
\cdots\, dx_{p+q-2}=0.
\end{eqnarray}
\end{thmm}

Rosi\'nski and Samorodnitsky \cite{RS} observed that multiple
Wiener--It\^o integrals are independent if and only if their squares
are uncorrelated,
%
\begin{equation}
\label{eqRS} I_p(f) \perp\!\!\!\!\perp I_q(g) \quad\iff\quad\operatorname{Cov}
\bigl(I_p(f)^2, I_q(g)^2 \bigr)
=0.
\end{equation}
This condition can be viewed as a generalization of the usual
covariance criterion for the independence of jointly Gaussian random
variables (the case of $p=q=1$).

In the seminal paper \cite{NP}, Nualart and Peccati discovered the
following surprising central limit theorem.

\begin{thmm}[(Nualart--Peccati)]\label{NPthmIntro}
Let $F_n=I_q(f_n)$, where $q\geq2$ is fixed and $f_n\in L^2(\R^q_+)$
are symmetric. Assume also that $E[F_n^2]= 1$ for all $n$.
Then convergence in distribution of $(F_n)$ to the standard normal law
is equivalent to convergence of the fourth moment. That is, as $n\to
\infty$,
%
\begin{equation}
\label{eqNP} F_n\stackrel{\mathrm{law}} {\to}N(0,1)
\quad\Longleftrightarrow\quad
E \bigl[F_n^4 \bigr]\to3.
\end{equation}
\end{thmm}
Shortly afterwards, Peccati and Tudor \cite{PT} established a
multidimensional extension of Theorem \ref{NPthmIntro}. Since the
publication of these two important papers, many improvements and
developments on this theme have been considered. In particular, Nourdin
and Peccati \cite{NPchi} extended Theorem \ref{NPthmIntro} to the case
when the limit of $F_n$'s is a centered gamma distributed random
variable. We refer the reader to \cite{ivangiobook} for
further information and details of the above results.

A heuristic argument linking Theorems~\ref{UZthmIntro} and~\ref
{NPthmIntro} was given by Rosi\'nski (\cite{R}, pages 3--4), while
addressing a question of Albert Shiryaev. Namely, let $F$ and $G$ be
two i.i.d. centered random variables with fourth moment and unit variance.
The link comes via a simple formula,
\[
\tfrac12\operatorname{Cov} \bigl((F+G)^2,(F-G)^2 \bigr)=E
\bigl[F^4 \bigr]-3,
\]
criterion \eqref{eqRS},
as well as the celebrated Bernstein theorem that asserts that $F$ and~$G$ are Gaussian
if and only if $F+G$ and $F-G$ are independent. A~rigorous argument to carry through this idea is based on a
characterization of the asymptotic independence of multiple Wiener--It\^
o integrals, which is much more difficult to handle than the plain
independence, and may also be of an independent interest. The
covariance between the squares of multiple Wiener--It\^o integrals
plays the pivotal role in this characterization.

At this point we should also mention an extension of \eqref{eqRS} to
the multivariate setting.
Let $I$ be a finite set and $(q_i)_{i\in I}$ be a sequence of
nonnegative integers.
Let $F_{i}=I_{q_i}(f_{i})$ be a multiple Wiener--It\^o integral of
order $q_i$, $i\in I$.
Consider a partition of $I$ into disjoint blocks $I_k$, so that $I=\bigcup_{k=1}^d I_k$, and the resulting random vectors $(F_{i})_{i\in I_k}$,
$k=1,\ldots,d$. Then
%
\begin{eqnarray}
\label{eqRS-m} \bigl\{ (F_{i})_{i\in I_k}\dvtx k \le d \bigr\}
\mbox{ are independent } \Leftrightarrow \operatorname{Cov} \bigl(F_{i}^2,F_{j}^2
\bigr)=0
\nonumber
\\[-8pt]
\\[-8pt]
\eqntext{\forall i, j \mbox{ from different blocks}.}
\end{eqnarray}

The proof of this criterion is similar to the proof of (1.2) in
\cite{RS}.

In this paper, in Theorem \ref{mainblock}, we establish an asymptotic
version of \eqref{eqRS-m} characterizing the asymptotic
moment-independence between blocks of multiple Wiener--It\^o integrals.
As a consequence of this result, we deduce the fourth moment theorem of
Nualart and Peccati \cite{NP} in Theorem \ref{nupec}, its
multidimensional extension due to Peccati and Tudor \cite{PT} in
Theorem \ref{pectud} and some neat estimates on the speed of
convergence in Theorem \ref{pectud-bd}. Furthermore, we obtain new
multidimensional extension of a theorem of Nourdin and Peccati~\cite
{NPchi} in Theorem \ref{tNouPec}, and give another new result on the
bivariate convergence of vectors consisting of multiple Wiener--It\^o
integrals in Theorem \ref{tbivariateblock}. Proposition \ref{BM-DMT}
applies Theorem~\ref{tbivariateblock} to establish the limit process
for functions of short and long-range dependent stationary Gaussian
time series in the spirit of the celebrated Breuer--Major \cite{BM} and
Dobrushin--Major--Taqqu \cite{DM,Taqqu} theorems.
In Theorem \ref{thmMOO} we establish the asymptotic moment-independence
for discrete non-Gaussian chaoses using some techniques of Mossel,
O'Donnel and Oleszkiewicz \cite{MOO}.

The paper is organized as follows. In Section~\ref{SPre} we list
some basic facts from Gaussian analysis and prove some lemmas needed in
the present work. In particular, we establish Lemma \ref{lCS}, which  is
a~version of the Cauchy--Schwarz inequality well suited to deal with
contractions of functions; see \eqref{v2}. It is used in the proof of
the main result, Theorem \ref{mainblock}. Section~\ref{Smain} is
devoted to the main results on the asymptotic independence. Section~\ref{Sapps}
gives some immediate consequences and related applications
of the main result. Section~\ref{sf-apps}
provides further applications to the study of short and long-range
dependent stochastic processes and multilinear random forms in
non-Gaussian random variables.

\section{Preliminaries}\label{SPre}

We will give here some basic elements of Gaussian analysis that are in
the foundations of the present work. The reader is referred to the
books~\cite{ivangiobook,Nbook} for further details
and ommited proofs.

Let $\EuFrak H$ be a real separable Hilbert space. For any $q\geq1$
let $\EuFrak H^{\otimes q}$ be the $q$th tensor product of $\EuFrak H$
and denote
by $\EuFrak H^{\odot q}$ the associated $q$th symmetric tensor product.
We write $X=\{X(h),h\in\EuFrak H\}$ to indicate
an isonormal Gaussian process over
$\EuFrak H$, defined on some probability space $(\Omega,\mathcal{F},P)$.
This means that $X$ is a centered Gaussian family, whose covariance is
given in terms of the
inner product of $\EuFrak H$ by $E [ X(h)X(g) ] =\langle
h,g\rangle_{\EuFrak H}$. We also assume that
$\mathcal{F}$ is
generated by $X$.

For every $q\geq1$, let $\mathcal{H}_{q}$ be the $q$th Wiener chaos of $X$,
that is, the closed linear subspace of $L^2(\Omega,\mathcal{F},P)$
generated by the random variables of the type $\{H_{q}(X(h)),h\in
\EuFrak H,\llVert
h\rrVert  _{\EuFrak H}=1\}$, where $H_{q}$ is the $q$th Hermite polynomial
defined as
%
\begin{equation}
\label{hermitepol} H_q(x) = (-1)^q e^{{x^2}/{2}}
\frac{d^q}{dx^q} \bigl( e^{-{x^2}/{2}} \bigr).
\end{equation}
We write by convention $\mathcal{H}_{0} = \mathbb{R}$. For
any $q\geq1$, the mapping
%
\begin{equation}
\label{mapping} I_{q} \bigl(h^{\otimes q} \bigr)=H_{q}
\bigl(X(h) \bigr)
\end{equation}
can be extended to a
linear isometry between  the symmetric tensor product $\EuFrak
H^{\odot q}$
equipped with the  modified norm $\sqrt{q!}\llVert  \cdot\rrVert
_{\EuFrak H^{\otimes q}}$ and the $q$th Wiener chaos $\mathcal{H}_{q}$.
For $q=0$ we write $I_{0}(c)=c$, $c\in\mathbb{R}$.

It is well known (Wiener chaos expansion) that $L^2(\Omega,\mathcal{F},P)$
can be decomposed into the infinite orthogonal sum of the spaces
$\mathcal{H}_{q}$. Therefore, any square integrable random variable
$F\in L^2(\Omega,\mathcal{F},P)$ admits the following chaotic expansion:
%
\begin{equation}
F=\sum_{q=0}^{\infty}I_{q}(f_{q}),
\label{E}
\end{equation}
where $f_{0}=E[F]$, and the $f_{q}\in\EuFrak H^{\odot q}$, $q\geq1$, are
uniquely determined by $F$. For every $q\geq0$ we denote by $J_{q}$ the
orthogonal projection operator on the $q$th Wiener chaos. In
particular, if
$F\in L^2(\Omega,\mathcal{F},P)$ is as in (\ref{E}), then
$J_{q}F=I_{q}(f_{q})$ for every $q\geq0$.

Let $\{e_{k}, k\geq1\}$ be a complete orthonormal system in $\EuFrak H$.
Given $f\in\EuFrak H^{\odot p}$ and $g\in\EuFrak H^{\odot q}$, for every
$r=0,\ldots,p\wedge q$, the \textit{contraction} of $f$ and $g$ of
order $r$
is the element of $\EuFrak H^{\otimes(p+q-2r)}$ defined by
%
\begin{equation}
f\otimes_{r}g=\sum_{i_{1},\ldots,i_{r}=1}^{\infty}
\langle f,e_{i_{1}}\otimes\cdots\otimes e_{i_{r}}
\rangle_{\EuFrak H^{\otimes
r}}\otimes\langle g,e_{i_{1}}\otimes\cdots\otimes
e_{i_{r}} \rangle_{\EuFrak H^{\otimes r}}. \label{v2}
\end{equation}
Notice that $f\otimes_{r}g$ is not necessarily symmetric: we denote its
symmetrization by $f\,\widetilde{\otimes}_{r}\,g\in\EuFrak H^{\odot(p+q-2r)}$.
Moreover, $f\otimes_{0}g=f\otimes g$ equals the tensor product of $f$ and
$g$ while, for $p=q$, $f\otimes_{q}g=\langle f,g\rangle_{\EuFrak
H^{\otimes q}}$.
In the particular case where $\EuFrak H=L^2(A,\mathcal{A},\mu)$, where
$(A,\mathcal{A})$ is a measurable space and $\mu$ is a $\sigma
$-finite and
nonatomic measure, one has that $\EuFrak H^{\odot q}=L_{s}^{2}(A^{q},
\mathcal{A}^{\otimes q},\mu^{\otimes q})$ is the space of symmetric and
square integrable functions on $A^{q}$. Moreover, for every $f\in
\EuFrak H^{\odot q}$, $I_{q}(f)$ coincides with the $q$-tuple
Wiener--It\^{o} integral of $f$.
In this case, (\ref{v2}) can be written as
\begin{eqnarray*}
&&(f\otimes_{r}g) (t_1,\ldots,t_{p+q-2r})\\
&&\qquad=\int
_{A^{r}}f(t_{1},\ldots,t_{p-r},s_{1},
\ldots,s_{r})
\\
&&\hspace*{14pt}\qquad\quad{}\times g(t_{p-r+1},\ldots,t_{p+q-2r},s_{1},
\ldots,s_{r})\,d\mu (s_{1})\cdots d\mu(s_{r}).
\end{eqnarray*}
We have
%
\begin{equation}
\label{fubini} \|f\otimes_r g\|^2 = \langle f
\otimes_{p-r} f, g \otimes_{q-r} g\rangle \qquad\mbox{for } r = 0,
\ldots, p\wedge q,
\end{equation}
where $\langle\cdot\rangle$ ($\|\cdot\|$, resp.) stands for inner
product (the norm, resp.) in an appropriate tensor product
space $\EuFrak H^{\otimes s}$.
Also, the following \textit{multiplication formula} holds: if $f\in\EuFrak
H^{\odot p}$ and $g\in\EuFrak
H^{\odot q}$, then
%
\begin{equation}
\label{multiplication} I_p(f) I_q(g) = \sum
_{r=0}^{p \wedge q} r! \pmatrix{p \cr r} \pmatrix{ q \cr r}
I_{p+q-2r} (f\,\widetilde{\otimes}_{r}\,g),
\end{equation}
where $f\,\widetilde{\otimes}_{r}\,g$ denotes the symmetrization of $f
\otimes_{r}g$.

We conclude these preliminaries with three useful lemmas that will be
needed throughout the sequel.
%
\begin{lemma}\label{usefullemma1}
\begin{longlist}[(ii)]
\item[(i)] Multiple Wiener--It\^o integral has all moments satisfying
the following hypercontractivity-type inequality:
%
\begin{equation}
\label{HC} \bigl[ E \bigl|I_p(f)\bigr|^r \bigr]^{1/r}
\le(r-1)^{p/2} \bigl[ E \bigl|I_p(f)\bigr|^2
\bigr]^{1/2},\qquad r\geq2.
\end{equation}
\item[(ii)] If a sequence of distributions of $\{I_p(f_n)\}_{n\ge1}$
is tight, then
%
\begin{equation}
\label{bounded} \sup_n E\bigl |I_p(f_n)\bigr|^r
< \infty\qquad \mbox{for every } r>0.
\end{equation}
\end{longlist}
\end{lemma}
\begin{pf}
(i) Inequality (\ref{HC}) is well known and corresponds,
for example, to~\cite{ivangiobook}, Corollary 2.8.14.

(ii) Combining (\ref{HC}) for $r=4$ with Paley's inequality, we get for
every $\theta\in(0,1)$,
%
\begin{equation}
\label{Paley} P \bigl( \bigl|I_p(f)\bigr|^2 > \theta E
\bigl|I_p(f)\bigr|^2 \bigr) \ge(1-\theta)^2
\frac
{ (E |I_p(f)|^2  )^2}{E |I_p(f)|^4} \ge(1-\theta)^2 9^{-p}.
\end{equation}
By the assumption, there is an $M>0$ such that $P(|I_p(f_n)|^2 > M) <
9^{-p-1}$, $n \ge1$. By \eqref{Paley} with $\theta=2/3$ and all $n$,
we have
\[
P \bigl(\bigl|I_p(f_n)\bigr|^2 > M \bigr) <
9^{-p-1} \le P \bigl(\bigl |I_p(f_n)\bigr|^2 >
(2/3) E \bigr|I_p(f_n)\bigr|^2 \bigr).
\]
As a consequence, $E |I_p(f_n)|^2 \le(3/2) M$. Applying \eqref{HC} we
conclude \eqref{bounded}.
\end{pf}

\begin{lemme}\label{technical}
\begin{longlist}[(1)]
\item[(1)] Let $p,q\geq1$, $f\in\EuFrak
H^{\odot p}$ and $g\in\EuFrak
H^{\odot q}$.
Then
%
\begin{equation}
\label{sym-norm} \|f\,\widetilde{\otimes}\,g\|^2 = \frac{p!q!}{(p+q)!}
\sum_{r=0}^{p\wedge
q} \pmatrix{p
\cr
r}\pmatrix{q
\cr
r} \|f\otimes_r g\|^2.
\end{equation}
\item[(2)] Let $q\geq1$ and $f_1,f_2,f_3,f_4\in\HH^{\odot q}$. Then
%
\begin{eqnarray}
\label{azerty} \qquad(2q)!\langle f_1\, \widetilde{\otimes}\,
f_2, f_3\, \widetilde{\otimes}\, f_4 \rangle &=&
\sum_{r=1}^{q-1} q!^2\pmatrix{q
\cr
r}^2 \langle f_1\otimes_r
f_3,f_4\otimes _r f_2\rangle
\nonumber
\\[-8pt]
\\[-8pt]
\nonumber
&&\hspace*{16pt}{}+q!^2 \bigl(\langle f_1,f_3\rangle \langle
f_2,f_4\rangle+ \langle f_1,f_4
\rangle \langle f_2,f_3\rangle \bigr).
\end{eqnarray}
\item[(3)] Let $q\geq1$, $f\in\HH^{\odot(2q)}$ and $g\in\HH^{\odot q}$.
We have
%
\begin{eqnarray}
\label{azertiop} &&\langle f\, \widetilde{\otimes}_q \,f, g \,\widetilde{
\otimes}\, g\rangle
\nonumber
\\[-8pt]
\\[-8pt]
\nonumber
&&\qquad= \frac{2q!^2}{(2q)!}\langle f\otimes_q f,g\otimes g
\rangle +\frac{q!^2}{(2q)!}\sum_{r=1}^{q-1}
\pmatrix{q
\cr
r}^2\langle f\otimes_r g,g
\otimes_r f\rangle.
\end{eqnarray}
\end{longlist}
\end{lemme}

\begin{pf} Without loss of generality, we suppose throughout the
proof that $\HH$ is equal to $L^2(A,\mathcal{A},\mu)$,
where $(A,\mathcal{A})$ is a measurable space, and $\mu$ is a $\sigma
$-finite measure without
atoms.

(1) Let $\sigma$ be a permutation of $\{1,\ldots,p+q\}$
(this fact is written in symbols as $\sigma\in\mathfrak{S}_{p+q}$).
If $r\in\{0,\ldots,p\wedge q\}$ denotes the cardinality of $\{1,\ldots
,p\}\cap\{\sigma(p+1),\ldots,\sigma(p+q)\}$,
then it is readily checked
that $r$ is also the cardinality of
$\{p+1,\ldots,p+q\}\cap\{\sigma(1),\ldots,\sigma(p)\}$ and that
%
\begin{eqnarray}
\label{ctr} &&\int_{A^{p+q}}f(t_1,
\ldots,t_p)g(t_{p+1},\ldots,t_{p+q})\nonumber\\
&&\hspace*{27pt}{}\times f(t_{\sigma
(1)},
\ldots,t_{\sigma(p)}) g(t_{\sigma(p+1)},\ldots,t_{\sigma(p+q)})\,d
\mu(t_1)\ldots d\mu (t_{p+q})
\nonumber
\\[-8pt]
\\[-8pt]
\nonumber
&&\qquad=\int_{A^{p+q-2r}}(f\otimes_r g) (x_1,
\ldots,x_{p+q-2r})^2 \,d\mu(x_1)\,\cdots\, d
\mu(x_{p+q-2r})\\
&&\qquad = \|f\otimes_r g\|^2.\nonumber
\end{eqnarray}
Moreover, for any fixed $r\in\{0,\ldots,p\wedge q\}$, there are
$p!{p\choose r} q!{q\choose r}$
permutations $\sigma\in\mathfrak{S}_{p+q}$ such that
$\{1,\ldots,p\}\cap\{\sigma(p+1),\ldots,\sigma(p+q)\}=r$.
[Indeed, such a permutation is completely determined by the choice of:
(a) $r$ distinct
elements $y_1,\ldots,y_r$ of $\{p+1,\ldots,p+q\}$; (b) $p-r$ distinct
elements $y_{r+1},\ldots,y_{p}$
of $\{1,\ldots,p\}$; (c) a bijection between $\{1,\ldots,p\}$ and $\{
y_1,\ldots,y_p\}$;
(d) a bijection between $\{p+1,\ldots,p+q\}$ and $\{1,\ldots,p+q\}
\setminus\{y_{1},\ldots,y_p\}$.]
Now, observe that the symmetrization of $f\otimes g$ is given by
\[
f\,\widetilde{\otimes}\, g(t_1,\ldots,t_{p+q}) =
\frac{1}{(p+q)!} \sum_{\sigma\in\mathfrak{S}_{p+q}} f(t_{\sigma(1)},
\ldots,t_{\sigma(p)}) g(t_{\sigma(p+1)},\ldots,t_{\sigma(p+q)}).
\]
Therefore, using (\ref{ctr}), we can write
\begin{eqnarray*}
\|f\,\widetilde{\otimes}\, g\|^2&=&\langle f\otimes g,f\,\widetilde{
\otimes}\, g\rangle \\
&=& \frac{1}{(p+q)!} \sum_{\sigma\in\mathfrak{S}_{p+q}} \int
_{A^{p+q}} f(t_{1},\ldots,t_{p})g(t_{p+1},
\ldots,t_{p+q})
\\
&&\hspace*{98pt}{} \times f(t_{\sigma(1)},\ldots,t_{\sigma(p)})\\
&&\hspace*{98pt}{} \times g(t_{\sigma(p+1)},
\ldots,t_{\sigma(p+q)}) \,d\mu(t_1)\cdots d\mu(t_{p+q})
\\
&=&\frac{1}{(p+q)!}\sum_{r=0}^{p\wedge q} \|f
\otimes_r g\|^2 \operatorname{Card} \bigl\{\sigma\in
\mathfrak{S}_{p+q}\dvtx \{1,\ldots,p\}\\
&&\hspace*{133pt}{}\cap \bigl\{\sigma(p+1),\ldots,
\sigma(p+q) \bigr\}=r \bigr\},
\end{eqnarray*}
and (\ref{sym-norm}) follows.

(2) We proceed analogously.
Indeed, we have
\begin{eqnarray*}
&&\langle f_1 \,\widetilde{\otimes}\, f_2,f_3
\,\widetilde{\otimes}\, f_4\rangle\\
&&\qquad= \langle f_1\otimes
f_2,f_3 \,\widetilde{\otimes}\, f_4\rangle
\\
&&\qquad= \frac{1}{(2q)!} \sum_{\sigma\in\mathfrak{S}_{2q}} \int
_{A^{2q}} f_1(t_{1},
\ldots,t_{q})f_2(t_{q+1},\ldots,t_{2q})
\\
&&\hspace*{106pt}{} \times f_3(t_{\sigma(1)},\ldots,t_{\sigma(q)})\\
&&\hspace*{106pt}{}\times f_4(t_{\sigma(q+1)},
\ldots,t_{\sigma(2q)}) \,d\mu(t_1)\cdots d\mu(t_{2q})
\\
&&\qquad=\frac{1}{(2q)!}\sum_{r=0}^q \langle
f_1\otimes_r f_3,f_4
\otimes_r f_2\rangle\\
&&\hspace*{73pt}{}\times \operatorname{Card} \bigl\{\sigma\in
\mathfrak{S}_{2q}\dvtx \bigl\{\sigma(1),\ldots,\sigma(q) \bigr\} \cap
\{1, \ldots,q\}=r \bigr\},
\end{eqnarray*}
from which we deduce (\ref{azerty}).

(3) We have
\begin{eqnarray*}
&&(g \,\widetilde{\otimes}\, g) (t_1,\ldots,t_{2q})\\
&&\qquad=
\frac{1}{(2q)!}\sum_{\sigma\in\mathfrak{S}_{2q}} g(t_{\sigma(1)},
\ldots,t_{\sigma(q)}) g(t_{\sigma(q+1)},\ldots,t_{\sigma(2q)})
\\
&&\qquad=\frac{1}{(2q)!}\sum_{r=0}^{q}
\mathop{\sum_{\sigma\in\mathfrak{S}_{2q}}}_{
\{\sigma(1),\ldots,\sigma(q)\}\cap\{1,\ldots,q\}=r
}
g(t_{\sigma(1)},\ldots,t_{\sigma(q)}) g(t_{\sigma(q+1)},
\ldots,t_{\sigma(2q)})
\end{eqnarray*}
and
\begin{eqnarray*}
&&(f\otimes_q f) (t_1,\ldots,t_{2q})= \int
_{A^{q}}f(t_1,\ldots,t_q,x_1,
\ldots,x_q)\\
&&\hspace*{124pt}{}\times f(x_1,\ldots,x_q,t_{q+1},
\ldots,t_{2q})\,d\mu(x_1)\cdots d\mu(x_q),
\end{eqnarray*}
so that
\begin{eqnarray*}
&&\langle f \,\widetilde{\otimes}_q\, f,g \,\widetilde{\otimes}\, g\rangle \\
&&\qquad=
\langle f\otimes_q f,g \,\widetilde{\otimes}\, g\rangle
\\
&&\qquad=\frac{1}{(2q)!}\sum_{r=0}^q \langle
f\otimes_r g,g\otimes_r f\rangle \\
&&\hspace*{75pt}{}\times \operatorname{Card} \bigl\{\sigma
\in\mathfrak{S}_{2q}\dvtx \bigl\{\sigma(1),\ldots,\sigma(q) \bigr\}\cap
\{1,\ldots,q\}=r \bigr\}
\\
&&\qquad= \frac{1}{(2q)!} \sum_{r=0}^q
\pmatrix{q
\cr
r}^2 q!^2 \langle f\otimes_r
g,g \otimes_r f\rangle
\\
&&\qquad= \frac{q!^2}{(2q)!}\langle f\otimes_q g,g\otimes_q f
\rangle+ \frac{q!^2}{(2q)!}\langle f\otimes g,g\otimes f\rangle \\
&&\qquad\quad{}+ \frac1{(2q)!}\sum
_{r=1}^{q-1}\pmatrix{q
\cr
r}^2q!^2\langle f\otimes _rg,g
\otimes_rf\rangle.
\end{eqnarray*}
Since $\langle f\otimes_q g,g\otimes_q f\rangle=\langle f\otimes
g,g\otimes f\rangle=\langle f\otimes_qf,g\otimes g\rangle$, the
desired conclusion~(\ref{azertiop}) follows.
\end{pf}

\begin{lemme}[(Generalized Cauchy--Schwarz inequality)]\label{lCS}
Assume that $\HH= L^2(A,\mathcal{A},\mu)$,
where $(A,\mathcal{A})$ is a measurable space equipped with a $\sigma
$-finite measure $\mu$.
For any integer $M\geq1$, put $[M] = \{1,\ldots,M\}$.
Also, for every element $\mathbf{z} = (z_1,\ldots,z_M)\in A^M$ and every
nonempty set $c\subset[M]$, let $\mathbf{z}_c$ denote the element of $A^{|c|}$
(where $|c|$ is the cardinality of $c$) obtained by deleting from $\mathbf{
z}$ the entries with index not contained in $c$.
(E.g., if $M=5$ and $c=\{1,3,5\}$, then $\mathbf{z}_c = (z_1,z_3,z_5)$.)
Let:
\begin{longlist}[{($\alpha$)}]
\item[{($\alpha$)}] $C,q\geq2$ be integers, and let $c_1,\ldots,c_q$
be nonempty subsets of $[C]$ such that each element of $[C]$ appears in
exactly two of the $c_i$'s (this implies that $\bigcup_{i} c_i = [C]$
and $\sum_i |c_i| = 2C$);

\item[{($\beta$)}] let $h_1,\ldots,h_q$ be functions such that $h_i
\in L^2(\mu^{|c_i|}):= L^2(A^{|c_i|},\mathcal{A}^{|c_i|},\mu^{|c_i|})$
for every $i=1,\ldots,q$ (in particular, each $h_i$ is a function of
$|c_i|$ variables).
\end{longlist}
Then
%
\begin{equation}
\label{genCS} \Biggl|\int_{A^C} \prod
_{i=1}^q h_i(\mathbf{z}_{c_i})
\mu^{C}(d\mathbf{ z}_{[C]}) \Biggr| \leq\prod
_{i=1}^q \|h_i\|_{L^2(\mu^{|c_i|})}.
\end{equation}
Moreover, if $c_0:= c_j \cap c_k \ne\varnothing$ for some $j \ne k$, then
%
\begin{equation}
\label{genCS1} \qquad\Biggl|\int_{A^C} \prod
_{i=1}^q h_i(\mathbf{z}_{c_i})
\mu^{C}(d\mathbf{ z}_{[C]}) \Biggl|\leq\|h_j
\otimes_{c_0} h_k\|_{L^2(\mu^{|c_j \triangle
c_k|})} \prod
_{i\ne j, k}^q \|h_i\|_{L^2(\mu^{|c_i|})},
\end{equation}
where
\[
h_j \otimes_{c_0} h_k(\mathbf{z}_{c_j \triangle c_k})
= \int_{A^{|c_0|}} h_j(\mathbf{z}_{c_j})
h_k(\mathbf{z}_{c_k}) \mu^{|c_0|}(d
\mathbf{z}_{c_0}).
\]
(Notice that $h_j \otimes_{c_0} h_k = h_j \otimes_{|c_0|} h_k$ when
$h_j$ and $h_k$ are symmetric.)
\end{lemme}

\begin{pf}
In the case $q=2$, \eqref{genCS} is just the Cauchy--Schwarz
inequality, and~\eqref{genCS1} is an equality. Assume that
\eqref{genCS}--\eqref{genCS1} hold for at most $q-1$ functions and proceed by
induction. Among the sets $c_1,\ldots,c_q$ at least two, say $c_j$ and~$c_k$, have nonempty intersections. Set $c_0:= c_j \cap c_k$, as
above. Since $c_0$ does not have common elements with $c_i$ for all $i
\ne j, k$, by Fubini's theorem,
%
\begin{eqnarray}
\label{genCS2} &&\int_{A^C} \prod
_{i=1}^q h_i(\mathbf{z}_{c_i})
\mu^{C}(d \mathbf{z}_{[C]})
\nonumber
\\[-8pt]
\\[-8pt]
\nonumber
&&\qquad = \int_{A^{C-|c_0|}}
h_j \otimes_{c_0} h_k(\mathbf{z}_{c_j \triangle c_k})
\prod_{i\ne j, k}^q h_i(
\mathbf{z}_{c_i}) \mu^{C-|c_0|}(d\mathbf{z}_{[C]
\setminus c_0}).
\end{eqnarray}
Observe that every element of $[C] \setminus c_0$ belongs to exactly
two of the $q-1$ sets: $c_j \triangle c_k$, $c_i$, $i\ne j,k$.
Therefore, by the induction assumption, \eqref{genCS} implies \eqref
{genCS1}, provided $c_j \triangle c_k \ne\varnothing$. When $c_j = c_k$,
we have $h_j \otimes_{c_0} h_k= \langle h_j, h_k \rangle$, and \eqref
{genCS1} follows from \eqref{genCS} applied to the product of $q-2$
functions in \eqref{genCS2}. This proves~\eqref{genCS1}, which in turn
yields \eqref{genCS} by the Cauchy--Schwarz inequality. The proof is complete.
\end{pf}

\section{The main results}\label{Smain}

The following theorem characterizes moment-inde\-pendence of limits of
multiple Wiener--It\^o integrals.

\begin{thmm}\label{main}
Let $d\geq2$, and let $q_1,\ldots,q_d$ be positive integers.
Consider vectors
\[
(F_{1,n}, \ldots, F_{d,n}) = \bigl(I_{q_1}(f_{1,n}),
\ldots, I_{q_d}(f_{d,n}) \bigr),\qquad  n\geq1,
\]
with $f_{i,n}\in\HH^{\odot q_i}$. Assume that for some random vector
$(U_1,\ldots,U_d)$,
%
\begin{equation}
\label{joint} (F_{1,n}, \ldots, F_{d,n}) \stackrel{\mathrm{law}} {
\to} (U_1,\ldots,U_d) \qquad\mbox{as } n \to\infty.
\end{equation}
Then $U_i$'s admit moments of all orders and the following three
conditions are equivalent:
\begin{longlist}[($\alpha$)]
\item[($\alpha$)] $U_1,\ldots,U_d$ are moment-independent, that is,
$E[U_1^{k_1}\cdots  U_d^{k_d}] = E[U_1^{k_1}]\cdots\break  E[U_d^{k_d}]$   for
all $k_1,\ldots,k_d \in\mathbb{N}$;
\item[($\beta$)] $\lim_{n\to\infty}\operatorname{Cov}(F_{i,n}^2,F_{j,n}^2) = 0$
for all $i\neq j$;
\item[($\gamma$)] $\lim_{n\to\infty}\|f_{i,n}\otimes_r f_{j,n}\| = 0$
for all $i\neq j$ and all $r=1,\ldots,q_i\wedge q_j$.
\end{longlist}
Moreover, if the distribution of each $U_i$ is determined by its
moments, then (a) is equivalent to that:
\begin{longlist}
\item[($\delta$)] $U_1,\ldots,U_d$ are independent.
\end{longlist}
\end{thmm}

\begin{rems}\label{rm1}
\begin{longlist}[(1)]
\item[(1)]
Theorem \ref{main} raises the question of whether the
moment-independence implies the usual independence under weaker
conditions than the determinacy of the marginals. (Recall that a random
variable having all moments is said to be determinate if any other
random variable with the same moments has the same distribution.) The
answer is negative in general; see \cite{BS}, Theorem 5.
%
%
%

\item[(2)] Assume that $d=2$ (for simplicity). In this case, ($\gamma$) becomes
$\|f_{1,n}\otimes_r f_{2,n}\|\to0$ for all $r=1,\ldots,q_1\wedge q_2$.
In view of Theorem \ref{UZthmIntro} of \"Ust\"unel and Zakai, one may
expect that ($\gamma$) could be replaced by a weaker condition ($\gamma
'$): $\|f_{1,n}\otimes_1 f_{2,n}\|\to0$.

However, the latter is false. To see it, consider a sequence $f_n \in\HH
^{\odot2}$ such that $\|f_n\|^2=\frac{1}{2}$ and $\|f_n\otimes_1 f_n\|
\to0$. By Theorem \ref{nupec} below, $F_n:= I_2(f_n)\stackrel{\mathrm{law}}{\to} U\sim N(0,1)$. Putting $f_{1,n}=f_{2,n}=f_n$, we observe that
($\gamma'$) holds, but ($\alpha$) does not, as
$(I_2(f_{1,n}),I_2(f_{2,n})) \stackrel{\mathrm{law}}{\to} (U,U)$.

\item[(3)]
Taking into account that assumptions $(\gamma)$ and $(\delta)$ of
Theorem \ref{nupec}
are equivalent, it is natural to wonder whether assumption ($\gamma$)
of Theorem \ref{main} is equivalent to its symmetrized version,
\[
\lim_{n\to\infty}\|f_{i,n}\, \widetilde{
\otimes}_r\, f_{j,n}\| = 0\qquad\mbox{for all $i\neq j$ and all
$r=1,\ldots,q_i\wedge q_j$}.
\]
The answer is negative in general, as is shown by the following counterexample.
Let $f_1,f_2\dvtx [0,1]^2\to\R$ be symmetric functions given by
\[
f_1(s,t)=\cases{ %
-1,& \quad $s,t\in[0,1/2],$
\vspace*{2pt}\cr
1,& \quad$\mbox{elsewhere}$}\quad \mbox{and}\quad
f_2(s,t)=\cases{ %
-1,&\quad  $s,t\in(1/2,1],$
\vspace*{2pt}\cr
1,& \quad$\mbox{elsewhere}.$}
\]
Then $\langle f_1,f_2\rangle=0$ and
\[
(f_1\otimes_1 f_2) (s,t)=\cases{
-1,& \quad $\mbox{if $s\in[0,1/2]$ and $t\in(1/2,1]$},$
\vspace*{2pt}\cr
1,&\quad $\mbox{if $t\in[0,1/2]$ and $s\in(1/2,1]$},$
\vspace*{2pt}\cr
0,&\quad $\mbox{elsewhere},$}
\]
so that $f_1\, \widetilde{\otimes}_1\, f_2\equiv0$
and $\|f_1\otimes_1 f_2\|=\sqrt{2}$.

\item[(4)]
The condition of moment-independence, ($\alpha$) of Theorem \ref{main},
can also be stated in terms of cumulants. Recall that the joint
cumulant of random variables $X_1,\ldots,X_m$ is defined by
\[
\kappa(X_1,\ldots,X_m) = (-i)^{m}
\frac{\partial^m}{\partial t_1 \cdots
\partial t_m} \log E \bigl[ {e^{i(t_1X_1+\cdots+t_mX_m)} \bigr]}_{ | t_1=0,\ldots
, t_m=0},
\]
provided $E|X_1\cdots X_m|< \infty$. When all $X_i$'s are equal to $X$, then
$\kappa(X,\ldots,\break  X)= \kappa_m(X)$, the usual $m$th cumulant of $X$; see
\cite{LS}. Then Theorem \ref{main}($\alpha$) is equivalent to
\begin{longlist}
\item[($\alpha'$)] \textit{for all integers $1\le j_1<\cdots<j_k \le d$, $k
\ge2$, and $m_1,\ldots,m_k\ge1$}
%
\begin{equation}
\label{eqa} \kappa(\underbrace{U_{j_1},\ldots,U_{j_1}}_{m_1},
\ldots, \underbrace {U_{j_k},\ldots,U_{j_k}}_{m_k}) =
0.
\end{equation}
\end{longlist}
\end{longlist}
\end{rems}

Theorem \ref{main} was proved in the first version of this paper \cite
{NR1}. Our proof of the crucial implication $(\gamma) \Rightarrow
(\alpha)$ involved tedious combinatorial considerations. We are
thankful to an anonymous referee who suggested a shorter and more
transparent line of proof using Malliavin calculus. It significantly
reduced the amount of combinatorial arguments of the original version
but requires some basic facts from Malliavin calculus.
We incorporated the referee's suggestions and approach into the proof
of a more general Theorem \ref{mainblock}. Even though Theorem \ref
{main} becomes a special case of Theorem \ref{mainblock} (see Corollary
\ref{marginalblock}), we keep its original statement for a convenient reference.

\begin{defi}
For each $n\geq1$, let $F_n=(F_{i,n})_{i\in I}$ be a family of
real-valued random variables indexed by a finite set $I$.
Consider a partition of $I$ into disjoint blocks $I_k$, so that $I=\bigcup_{k=1}^d I_k$.
We say that vectors $(F_{i,n})_{i\in I_k}$, $k=1,\ldots,d$ are
asymptotically moment-independent if each $F_{i,n}$ admits
moments of all orders and for any sequence $(\ell_i)_{i\in I}$ of
nonnegative integers,
%
\begin{equation}
\label{blocks} \lim_{n\to\infty} \Biggl\{ E \biggl[\prod
_{i\in I} F_{i,n}^{\ell_i} \biggr] -\prod
_{k=1}^d E \biggl[\prod
_{i\in I_k} F_{i,n}^{\ell_i} \biggr] \Biggr\}=0.
\end{equation}
\end{defi}

The next theorem characterizes the asymptotic moment-independence
between blocks of multiple Wiener--It\^o integrals.

\begin{thmm}\label{mainblock}
Let $I$ be a finite set and $(q_i)_{i\in I}$ be a sequence of
nonnegative integers.
For each $n\geq1$, let $F_n=(F_{i,n})_{i\in I}$ be a family of
multiple Wiener--It\^o integrals, where $F_{i,n}=I_{q_i}(f_{i,n})$ with
$f_{i,n}\in\HH^{\odot q_i}$. Assume that for every $i\in I$,
%
\begin{equation}
\label{u-bound} \sup_{n} E \bigl[F_{i,n}^2
\bigr] < \infty.
\end{equation}
Given a partition of $I$ into disjoint blocks $I_k$, the following
conditions are equivalent:
\begin{longlist}[(\textbf{a})]
\item[(\textbf{a})] random vectors $(F_{i,n})_{i\in I_k}$, $k=1,\ldots,d$
are asymptotically moment-independent;
\item[(\textbf{b})] $\lim_{n\to\infty} \operatorname{Cov}(F_{i,n}^2,F_{j,n}^2)=0$
for every $i,j$ from different blocks;
\item[(\textbf{c})] $\lim_{n\to\infty} \|f_{i,n}\otimes_r f_{j,n}\|=0$ for
every $i,j$ from different blocks and $r=1,\ldots,q_i\wedge
q_j$.\vadjust{\goodbreak}
\end{longlist}
\end{thmm}

\begin{pf}
The implication $\mathbf{(a)\Rightarrow(b)}$ is obvious.

To show $\mathbf{(b)\Rightarrow(c)}$, fix $i,j$ belonging to different
blocks. By \eqref{multiplication} we have
\[
F_{i,n}F_{j,n}=\sum_{r=0}^{q_i\wedge q_j}
r!\pmatrix{q_i
\cr
r}\pmatrix{q_j
\cr
r}
I_{q_i+q_j-2r}(f_{i,n}\,\widetilde{\otimes}_r\,
f_{j,n}),
\]
which yields
\[
E \bigl[F_{i,n}^2F_{j,n}^2 \bigr] =
\sum_{r=0}^{q_i\wedge q_j} r!^2
\pmatrix{q_i
\cr
r}^2\pmatrix{q_j
\cr
r}^2(q_i+q_j-2r)! \|f_{i,n}\,
\widetilde{\otimes}_r\, f_{j,n}\|^2.
\]
Moreover,
\[
E \bigl[F_{i,n}^2 \bigr]E \bigl[F_{j,n}^2
\bigr] = q_i!q_j!\|f_{i,n}\|^2
\|f_{j,n}\|^2.
\]
Applying (\ref{sym-norm}) to the second equality below, we evaluate
$\operatorname{Cov}(F_{i,n}^2,F_{j,n}^2)$ as follows:
%
\begin{eqnarray}
\label{cov} \operatorname{Cov} \bigl(F_{i,n}^2,F_{j,n}^2
\bigr) &=& (q_i+q_j)!\|f_{i,n}\,\widetilde{
\otimes }\, f_{j,n}\|^2 - q_i!q_j!
\|f_{i,n}\|^2\|f_{j,n}\|^2
\nonumber\\
&&{} +\sum_{r=1}^{q_i\wedge q_j} r!^2
\pmatrix{q_i
\cr
r}^2\pmatrix{q_j
\cr
r}^2(q_i+q_j-2r)! \|f_{i,n}\,
\widetilde{\otimes}_r\, f_{j,n}\|^2
\\
\nonumber
&=& q_i!q_j! \sum_{r=1}^{q_i\wedge q_j}
\pmatrix{q_i
\cr
r}\pmatrix{q_j
\cr
r}\|
f_{i,n}\otimes_r f_{j,n}\|^2\\
&&{} +\sum
_{r=1}^{q_i\wedge q_j} r!^2
\pmatrix{q_i
\cr
r}^2\pmatrix{q_j
\cr
r}^2(q_i+q_j-2r)! \|f_{i,n}\,
\widetilde{\otimes}_r\, f_{j,n}\|^2
\nonumber
\\
\label{cov2}&\geq&\max_{r=1,\ldots,q_i\wedge q_j} \|f_{i,n}\otimes_r
f_{j,n}\|^2.
\end{eqnarray}
This bound yields the desired conclusion.

Now we will prove $\mathbf{(c)\Rightarrow(a)}$. We need to show \eqref
{blocks} for fixed $l_i$. Writing $F_{i,n}^{l_i}$ as $\underbrace
{F_{i,n}\times\cdots\times
F_{i,n}}_{l_i}$ and enlarging $I$ and $I_k$'s accordingly, we may and
do assume that all $l_i=1$. We will prove (\ref{blocks}) by induction
on $Q=\sum_{i\in I}q_i$. The formula holds when $Q=0$ or 1. Therefore,
take $Q \ge2$ and suppose that (\ref{blocks}) holds whenever $\sum_{i\in I}q_i \le Q-1$.

Fix $i_1\in I_1$ and set
\[
X_n=\prod_{i\in I_1 \setminus\{i_1\}} I_{q_i}(f_{i,n}),\qquad
Y_n=\prod_{j\in I \setminus I_1} I_{q_j}(f_{j,n}).
\]
Assume that $q_1\ge1$, otherwise the inductive step follows immediately.
Let $\delta$ denote the divergence operator in the sense of Malliavin
calculus, and let $D$ be the Malliavin derivative; see
\cite{Nbook}, Chapters~1.2--1.3. Using the duality relation
\cite{Nbook}, Definition 1.3.1(ii), and the product rule for the Malliavin derivative
\cite{DOP}, Theorem~3.4, we get
\begin{eqnarray*}
E \biggl[\prod_{i\in I} F_{i,n} \biggr] &=&E
\bigl[I_{q_{i_1}}(f_{i_1,n})X_nY_n \bigr]
=E \bigl[\delta \bigl(I_{q_{i_1}-1}(f_{i_1,n}) \bigr)X_nY_n
\bigr]
\\
&=&E \bigl[I_{q_{i_1}-1}(f_{i_1,n})\otimes_1
D(X_nY_n) \bigr]
\\
& =& E \bigl[ Y_n I_{q_{i_1}-1}(f_{i_1,n})
\otimes_1 DX_n \bigr] + E \bigl[ X_n
I_{q_{i_1}-1}(f_{i_1,n})\otimes_1 DY_n
\bigr]
\\
& =& A_n+B_n.
\end{eqnarray*}
First we consider $B_n$. Using the product rule for $DY_n$, we obtain
\begin{eqnarray*}
B_n &=& \sum_{j\in I \setminus I_1} E
\biggl[I_{q_{i_1}-1}(f_{i_1,n})\otimes_1
DF_{j,n} \prod_{i\in I \setminus\{
i_1, j\}} F_{i,n}
\biggr]
\\
&=& \sum_{j\in I \setminus I_1} q_j E \biggl[
I_{q_{i_1}-1}(f_{i_1,n})\otimes_1 I_{q_j-1}(f_{j,n})
\prod_{i\in I
\setminus\{i_1, j\}} F_{i,n} \biggr].
\end{eqnarray*}
By the multiplication formula (\ref{multiplication}) we have
\begin{eqnarray*}
&&I_{q_{i_1}-1}(f_{i_1,n})\otimes_1 I_{q_j-1}(f_{j,n})
\\
&&\qquad= \sum_{s=1}^{q_{i_1}\wedge q_j}(s-1)!
\pmatrix{q_{i_1}-1
\cr
s-1}\pmatrix{q_{j}-1
\cr
s-1}I_{q_{i_1}+q_j-2s}(f_{i_1,n}\,\widetilde{\otimes}_s\,
f_{j,n}).
\end{eqnarray*}
Since $i_1$ and $j$ belong to different blocks, condition (c) of the
theorem applied to the above expansion yields that
$I_{q_{i_1}-1}(f_{i_1,n})\otimes_1 I_{q_{j}-1}(f_{j,n})$ converges to
zero in~$L^2$. Combining this with \eqref{u-bound} and Lemma \ref
{usefullemma1}, we infer that $\lim_{n\to\infty}B_n= 0$.

Now we consider $A_n$. If $\operatorname{Card}(I_1)=1$, then $X_n=1$ by
convention and so $A_n=0$. Hence
\[
\lim_{n \to\infty} \Biggl\{E \biggl[\prod
_{i\in I} F_{i,n} \biggr] -E [F_{i_1,n} ] \prod
_{k=2}^d E \biggl[\prod
_{i\in I_k} F_{i,n} \biggr] \Biggr\} =\lim
_{n \to\infty} B_n=0.
\]
Therefore, we now assume that $\operatorname{Card}(I_1)\ge2$.
Write $A_n=E [Z_nY_n ]$, where
\begin{eqnarray*}
Z_n &=&I_{q_{i_1}-1}(f_{i_1,n})\otimes_1
DX_n
\\
&=&\sum_{i\in I_1 \setminus\{i_1\}} q_i I_{q_{i_1}-1}(f_{i_1,n})
\otimes_1 I_{q_{i}-1}(f_{i,n})\prod
_{j\in I_1 \setminus\{i_1, i\}} F_{j,n}
\\
&=&\sum_{i\in I_1 \setminus\{i_1\}}q_i \sum
_{s=1}^{q_{i_1}\wedge q_i}(s-1)!\pmatrix{q_{i_1}-1
\cr
s-1} \pmatrix{q_{i}-1
\cr
s-1}\\
&&\hspace*{44pt}\qquad{}\times I_{q_{i_1}+q_i-2s}(f_{i_1,n}
\,\widetilde{ \otimes}_s\, f_{i,n}) \prod
_{j\in I \setminus\{i_1, i\}} F_{j,n}.
\end{eqnarray*}
Thus $A_n$ is a linear combination of the terms
\[
E \biggl[ \biggl( I_{q_{i_1}+q_i-2s}(f_{i_1,n}\,\widetilde{
\otimes}_s\, f_{i,n}) \prod_{j\in I_1 \setminus\{i_1, i\}}
F_{j,n} \biggr) Y_n \biggr],
\]
where $i_1,i \in I_1$, $i_1\ne i$, $1 \le s \le q_{i_1} \wedge q_i$.
The term under expectation is a product of multiple integrals of orders
summing to $\sum_{j \in I} q_j - 2s$. Therefore, the induction
hypothesis applies provided
%
\begin{equation}
\label{toshow2} \lim_{n\to\infty} (f_{i_1,n}\,\widetilde{
\otimes}_s\, f_{i,n} )\otimes_r
f_{j,n}= 0
\end{equation}
for all $j\in I_k$ with $k\geq2$ and all $r=1,\ldots
,(q_{i_1}+q_i-2s)\wedge q_j$.

Suppose that \eqref{toshow2} holds. Then by the induction hypothesis,
\[
\lim_{n \to\infty} \bigl\{ A_n - E[Z_n]
E[Y_n] \bigr\} =0.
\]
Moreover,
\[
E[Z_n]=E \bigl[I_{q_{i_1}-1}(f_{i_1,n})
\otimes_1 DX_n \bigr]=E \bigl[I_{q_{i_1}}(f_{i_1,n})X_n
\bigr]= E \biggl[\prod_{i\in I_1} F_{i,n}
\biggr].
\]
Hence, by the induction hypothesis applied to $Y_n$ and the uniform
boundedness of all moments of $F_{i,n}$, we get
\[
\lim_{n\to\infty} \Biggl\{ E \biggl[\prod
_{i\in I} F_{i,n} \biggr] -\prod
_{k=1}^d E \biggl[\prod
_{i\in I_k} F_{i,n} \biggr] \Biggr\} = \lim
_{n\to\infty} \bigl\{A_n - E[Z_n]
E[Y_n] \bigr\} = 0.
\]

It remains to show \eqref{toshow2}. To this aim we will describe the
structure of the terms under limit \eqref{toshow2}.
Without loss of generality we may assume that $\HH=L^2(\mu
):=L^2(A,\mathcal{A},\mu)$, where $(A,\mathcal{A})$ is a measurable
space and $\mu$ is a $\sigma$-finite measure without atoms. Recall
notation of Lemma \ref{lCS}.
For every integer $M\geq1$, put $[M] = \{1,\ldots,M\}$. Also, for every
element $\mathbf{z} = (z_1,\ldots,z_M)\in A^M$ and every nonempty set
$c\subset[M]$,
we denote by $\mathbf{z}_c$ the element of $A^{|c|}$
(where $|c|$ is the cardinality of $c$) obtained by deleting from $\mathbf{
z}$ the entries with index not contained in $c$.
(E.g., if $M=5$ and $c=\{1,3,5\}$, then $\mathbf{z}_c = (z_1,z_3,z_5)$.)

Observe that $ (f_{i_1,n}\,\widetilde{\otimes}_s\, f_{i,n} )\otimes_r
f_{j,n}$ is a linear combination of functions $\psi(\mathbf{z}_{J_1})$,
$\mathbf{z} \in A^M$ obtained as follows. Set $M=q_{i_1}+q_i + q_j -s-r$
and $M_0= q_{i_1}+q_i-s$, so that $M > M_0 \ge2$.
Choose $b_1, b_2 \subset[M_0]$ such that $|b_1|=q_{i_1}$, $|b_2| =
q_i$ and $|b_1\cap b_2|=s$, and then choose $b_3 \subset[M]$ such that
$|b_3|=q_j$ and $|b_3 \cap(b_1 \cup b_2)|=r$. It follows that $b_1\cup
b_2 \cup b_3 = [M]$ and $b_1\cap b_2 \cap b_3 = \varnothing$. Therefore,
each element of $[M]$ belongs exactly to one or two $b_i$'s. Let
\[
J= \bigl\{j\in[M]\dvtx j \mbox{ belongs to two sets } b_i \bigr\}
\]
and put $J_1=[M] \setminus J$. Then $ (f_{i_1,n}\,\widetilde{\otimes
}_s\, f_{i,n} )\otimes_r f_{j,n}$ is a linear combination of functions
of the form
\[
\psi(\mathbf{z}_{J_1}) = \int_{A^J}
f_{i_1,n}(\mathbf{z}_{b_1}) f_{i,n}(\mathbf{
z}_{b_2}) f_{j,n}(\mathbf{z}_{b_3})
\mu^{|J|}(d\mathbf{z}_J),
\]
where the summation goes over all choices $b_1, b_2$ under the
constraint that the sets $b_1 \cap b_2$ and $b_3$ are fixed. This
constraint makes $J_1$ unique, $|J_1|=q_{i_1}+q_i + q_j -2s-2r$.

Let $c_i = b_i \cap J$, $i=1,2,3$ and notice that either $c_1\cap c_3
\ne\varnothing$ or $c_2\cap c_3 \ne\varnothing$ since $r \ge1$. Suppose
$c_0=c_1\cap c_3 \ne\varnothing$, the other case is identical. Applying
Lemma \ref{lCS} with $\mathbf{z}_{J_1}$ fixed, we get
\[
\bigl|\psi(\mathbf{z}_{J_1})\bigr|^2 \le\bigl|f_{i_1,n}
\otimes_{|c_0|}f_{j,n}(\mathbf{ z}_{b_1\triangle b_3})\bigr|^2
\int_{A^{|c_2|}} \bigl|f_{i,n}(\mathbf{z}_{b_2})\bigr|^2
\mu^{|c_2|}(d\mathbf{z}_{c_2}).
\]
Since $b_1 \triangle b_3$ and $b_3\setminus c_3$ make a
disjoint partition of $J_1$, and additional integration
with respect to $\mathbf{z}_{J_1}$ yields
\[
\|\psi\|_{L^2(\mu^{|J_1|})} \le\|f_{i_1,n}\otimes_{|c_0|}f_{j,n}
\| _{L^2(\mu^{|b_1\triangle b_3|})} \|f_{i,n}\|_{L^2(\mu^{|b_2|})} \to0
\]
as $n \to\infty$. This yields \eqref{toshow2} and completes the proof
of Theorem \ref{mainblock}.
\end{pf}

\begin{remark}\label{rm2}
Condition ($\mathbf{b}$) of Theorem \ref{mainblock} is equivalent to

\begin{longlist}
\item[(\textbf{b}$\bolds{'}$)] \textit{for every $1 \le k \ne l \le d$}
\[
\lim_{n\to\infty} \operatorname{Cov} \bigl(\bigl\Vert(F_{i,n})_
{i\in I_k}\bigr\Vert^2,
\bigl\Vert(F_{i,n})_{i\in I_{l}}\bigr\Vert^2 \bigr)=0,
\]
\end{longlist}
where $\norm{\cdot}$ denotes the Euclidean norms in $\R^{|I_k|}$ and $\R
^{|I_{l}|}$, respectively.
\begin{pf}
Indeed, condition (\textbf{b}) of Theorem \ref{mainblock} implies (\textbf{b}$\bolds{'}$), and the converse follows from
\[
\operatorname{Cov} \bigl(\bigl\Vert(F_{i,n})_{i\in I_k}\bigr\Vert^2,
\bigl\Vert(F_{i,n})_{i\in
I_{l}}\bigr\Vert^2 \bigr)= \sum
_{i\in I_k, j \in I_l} \operatorname {Cov} \bigl(F_{i,n}^2,F_{j,n}^2
\bigr) \ge\operatorname{Cov} \bigl(F_{i,n}^2,F_{j,n}^2
\bigr)
\]
as the squares of multiple Wiener--It\^o integrals are nonnegatively
correlated; cf.~\eqref{cov2}.
\end{pf}
\end{remark}

The following corollary is useful in deducing the joint convergence in
law from the convergence of marginals. It is stated for random vectors,
as is Theorem \ref{mainblock}, but it obviously applies in the setting
of Theorem \ref{main} when all vectors are one-dimensional.

\begin{cor}\label{marginalblock}
Under notation of Theorem \ref{mainblock}, let $(U_i)_{i \in I}$ be a
random vector such that:
\begin{longlist}[(iii)]
\item[(i)] $(F_{i,n})_{i\in I_k} \stackrel{\mathrm{law}}{\to} (U_i)_{i\in
I_k}$   as $n\to\infty$, for each $k$;

\item[(ii)] vectors $(U_i)_{i\in I_k}$, $k=1,\ldots,d$ are independent;

\item[(iii)] condition $(\mathbf{b})$ or $(\mathbf{c})$ of Theorem \ref
{mainblock} holds [equivalently, $(\beta)$ or $(\gamma)$ of Theorem \ref
{main} when all $I_k$ are singletons];

\item[(iv)] $\mathcal{L}(U_i)$ is determined by its moments for each $i
\in I$.
\end{longlist}
Then the joint convergence holds
\[
(F_{i,n})_{i\in I} \stackrel{\mathrm{law}} {\to}
(U_i)_{i\in I},\qquad n \to \infty.
\]
\end{cor}

\begin{pf}
By (i) the sequence $\{(F_{i,n})_{i\in I}\}_{n \ge1}$ is tight.
Let $(V_i)_{i \in I}$ be a random vector such that
\[
(F_{i, n_j})_{i \in I} \stackrel{\mathrm{law}} {\to}
(V_i)_{i \in I}
\]
as $n_j \to\infty$ along a subsequence. From Lemma \ref
{usefullemma1}(ii) we infer that condition~\eqref{u-bound} of Theorem
\ref{mainblock} is satisfied. It follows that each $V_i$ has all
moments and
$(V_i )_{i\in I_k}   \stackrel{\mathrm{law}}{=}   (U_i)_{i\in I_k}$ for
each $k$. By (iv), the laws of vectors $(U_i)_{i\in I}$ and
$(V_i)_{i\in I}$ are determined by their joint moments, respectively;
see \cite{P}, Theorem 3. Under assumption (iii), the vectors
$(F_{i,n})_{i\in I_k}$, $k=1,\ldots,d$ are asymptotically moment
independent. Hence,
for any sequence $(\ell_i)_{i\in I}$ of nonnegative integers,
\begin{eqnarray*}
E \biggl[\prod_{i\in I} V_{i}^{\ell_i}
\biggr] - E \biggl[\prod_{i\in I} U_{i}^{\ell_i}
\biggr] &=& E \biggl[\prod_{i\in I} V_{i}^{\ell_i}
\biggr] - \prod_{k=1}^d E \biggl[\prod
_{i\in I_k} U_{i}^{\ell_i} \biggr]
\\
&=& \lim_{n_j\to\infty} \Biggl\{ E \biggl[\prod
_{i\in I} F_{i,n_j}^{\ell_i} \biggr] -\prod
_{k=1}^d E \biggl[\prod
_{i\in I_k} F_{i,n_j}^{\ell_i} \biggr] \Biggr\}=0.
\end{eqnarray*}
Thus $(V_i )_{i\in I}   \stackrel{\mathrm{law}}{=}   (U_i)_{i\in I}$.
\end{pf}

\section{Applications} \label{Sapps}

\subsection{The fourth moment theorem of Nualart--Peccati}

We can give a short proof of the difficult and surprising part
implication $(\beta) \Rightarrow(\alpha)$ of the fourth moment theorem
of Nualart and Peccati \cite{NP}, that we restate here for a convenience.

\begin{thmm}[(Nualart--Peccati)]\label{nupec}
Let $(F_n)$ be a sequence of the form $F_n=I_q(f_n)$, where $q\geq2$
is fixed and $f_n\in\HH^{\odot q}$.
Assume moreover that $E[F_n^2]=q!\|f_n\|^2=1$ for all $n$.
Then, as $n\to\infty$, the following four conditions are equivalent:
\begin{longlist}[($\alpha$)]
\item[($\alpha$)] $F_n\stackrel{\mathrm{law}}{\to}N(0,1)$;
\item[($\beta$)] $E[F_n^4]\to3$;\vadjust{\goodbreak}
\item[($\gamma$)] $\|f_{n}\otimes_r f_{n}\|\to0$ for all $r=1,\ldots
,q-1$;
\item[($\delta$)] $\|f_{n}\, \widetilde{\otimes}_r \, f_{n}\|\to0$ for
all $r=1,\ldots,q-1$.
\end{longlist}
\end{thmm}

\begin{pf*}{Proof of $(\beta) \Rightarrow(\alpha)$} Assume $(\beta)$.
Since the sequence $(F_n)$ is bounded in $L^2(\Omega)$ by the
assumption, it is relatively compact in law. Without loss of generality
we may assume that $F_n\stackrel{\mathrm{law}}{\to} Y$ and need to show that
$Y \sim N(0,1)$.
Let $G_n$ be an independent copy of $F_n$ of the form $G_n=I_q(g_n)$
with $f_n\otimes_1 g_n=0$. This can easily be done by extending the
underlying isonormal process to the direct sum $\EuFrak H \oplus
\EuFrak H$.
We then have
\[
\bigl(I_{q} (f_{n}+g_{n} ),I_q
(f_{n}-g_{n} ) \bigr) = (F_{n}+G_{n},F_{n}-G_{n}
) \stackrel{\mathrm{law}} {\to} (Y+Z,Y-Z )
\]
as $n\to\infty$, where $Z$ stands for an independent copy of $Y$.
Since
\[
\tfrac{1}{2} \operatorname{Cov} \bigl[(F_{n}+G_{n})^2,
(F_{n}-G_{n})^2 \bigr] = E
\bigl[F_{n}^4 \bigr] - 3 \to0,
\]
$Y+Z$ and $Y-Z$ are moment-independent. (If they were independent, the
classical Bernstein theorem would complete the proof.) However, in our
case condition ($\alpha'$) in \eqref{eqa} says that
\[
\kappa(\underbrace{Y+Z,\ldots, Y+Z}_{m_1}, \underbrace{Y-Z,\ldots,
Y-Z}_{m_2}) = 0\qquad \mbox{for all } m_1,m_2 \ge1.
\]
Taking $n\ge3$ we get
\begin{eqnarray*}
0 &= &\kappa(\underbrace{Y+Z,\ldots, Y+Z}_{n-2}, Y-Z, Y-Z)
\\
& =& \kappa(\underbrace{Y,\ldots, Y}_{n}) + \kappa(\underbrace{Z,
\ldots, Z}_{n}) = 2 \kappa_n(Y),
\end{eqnarray*}
where we used the multilinearity of $\kappa$ and the fact that $Y$ and
$Z$ are i.i.d. Since $\kappa_1(Y)=0$, $\kappa_2(Y)=1$ and $\kappa
_n(Y)=0$ for $n \ge3$, we infer that $Y \sim N(0,1)$.
\end{pf*}

\subsection{Generalizing a result of Peccati and Tudor}

Applying our approach, one can add a further equivalent condition to a
result of Peccati and Tudor \cite{PT}. As such, Theorem \ref{pectud}
turns out to be the exact multivariate equivalent of Theorem~\ref{nupec}.

\begin{thmm}[(Peccati--Tudor)]\label{pectud}
Let $d\geq2$, and let $q_1,\ldots, q_d$ be positive integers.
Consider vectors
\[
F_n=(F_{1,n},\ldots,F_{d,n}) =
\bigl(I_{q_1}(f_{1,n}),\ldots, I_{q_d}(f_{d,n})
\bigr),\qquad n\geq1,
\]
with $f_{i,n}\in\HH^{\odot q_i}$. Assume that, for $i,j=1,\ldots,d$, as
$n\to\infty$,
%
\begin{equation}
\label{kron} \operatorname{Cov} (F_{i,n}, F_{j,n} ) \to
\sigma_{ij}.
\end{equation}
Let $N$ be a centered Gaussian random vector with the covariance matrix
$\Sigma=(\sigma_{ij})_{1\leq i,j\leq d}$.
Then the following two conditions are equivalent ($n\to\infty$):
\begin{longlist}[(ii)]
\item[(i)] $F_n\stackrel{\mathrm{law}}{\to}N$;
\item[(ii)]
$E [\|F_n\|^4 ]\to E [\|N\|^4 ]$;
\end{longlist}
where $\|\cdot\|$ denotes the Euclidean norm in $\R^d$.
\end{thmm}

\begin{pf}
Only $\mathrm{(ii)}\Rightarrow\mathrm{(i)}$ has to be shown. Assume
(ii).
As in the proof of Theorem \ref{nupec}, we may assume that $F_n\stackrel
{\mathrm{law}}{\to} Y$ and must show that $Y \sim N_d(0, \Sigma)$. Let
$G_n=(G_{1,n},\ldots,G_{d,n})$ be an independent copy of $F_n$ of the
form $ (I_{q_1}(g_{1,n}),\ldots, I_{q_d}(g_{d,n}) )$. Observe that
\begin{eqnarray*}
&&\frac{1}{2}\operatorname{Cov} \bigl(\norm{F_{n}+G_{n}}^2
, \norm {F_{n}-G_{n}}^2 \bigr)
\\
&&\qquad = E \bigl[\norm{F_{n}}^4 \bigr] - \bigl(E \bigl[
\norm{F_{n}}^2 \bigr] \bigr)^2 - 2 \sum
_{i,j=1}^d \operatorname{Cov}(F_{i,n},
F_{j,n})^2.
\end{eqnarray*}
Using this identity for $N$ and $N'$ in place of $F_n$ and $G_n$, where
$N'$ is an independent copy of $N$, we get
%
\begin{equation}
\label{identity} E \bigl[\norm{N}^4 \bigr] = \sum
_{i,j=1}^d \bigl( \sigma_{ii}
\sigma_{jj} + 2 \sigma _{ij}^2 \bigr).
\end{equation}
Hence
\begin{eqnarray*}
&&\frac{1}{2}\operatorname{Cov} \bigl(  \norm{F_{n} +G_{n}}^2
, \norm {F_{n}-G_{n}}^2 \bigr) \\
&&\qquad= E \bigl[
\norm{F_{n}}^4 \bigr] -E \bigl[\norm{N}^4
\bigr]
\\
&&\qquad\quad{} + \sum_{i,j=1}^d \bigl[
\sigma_{ii} \sigma_{jj} + 2 \sigma_{ij}^2
- \operatorname{Var}(F_{i,n})\operatorname{Var}(F_{j,n}) - 2
\operatorname{Cov}(F_{i,n}, F_{j,n})^2 \bigr] \to0.
\end{eqnarray*}
By Remark \ref{rm2}, $F_n+G_n$ and $F_n-G_n$ are asymptotically
moment-independent. Since one-dimensional projections of $F_n+G_n$ and
$F_n-G_n$ are also asymptotically moment-independent, we can proceed by
cumulants as above to determine the normality of $Y$.
\end{pf}

The following result associates neat estimates to Theorem \ref{pectud}.

\begin{thmm}\label{pectud-bd}
Consider a vector
\[
F=(F_1,\ldots,F_d)= \bigl(I_{q_1}(f_1),
\ldots,I_{q_d}(f_d) \bigr)
\]
with $f_i\in\HH^{\odot q_i}$, and let $\Sigma=(\sigma_{ij})_{1\leq
i,j\leq d}$ be the covariance matrix of $F$, $\sigma_{ij}=E[F_iF_j]$.
Let $N$ be the associated Gaussian random vector, $N\sim N_d(0,\Sigma)$.
\begin{longlist}[(1)]
\item[(1)] Assume that $\Sigma$ is invertible.
Then, for any Lipschitz function $h\dvtx \R^d\to\R$ we have
\[
\bigl|E \bigl[h(F) \bigr] - E \bigl[h(N) \bigr] \bigr| \leq\sqrt{d} \|\Sigma
\|_{\mathrm{op}}^{1/2}\|\Sigma^{-1}\|_{\mathrm{op}}\|h
\|_{\mathrm{Lip}} \sqrt{E\|F\|^4-E\|N\|^4},\vadjust{\goodbreak}
\]
where $\|\cdot\|_{\mathrm{op}}$ denotes the operator norm of a matrix and $\|
h\|_{\mathrm{Lip}}=\break  \sup_{x,y\in\R^d}\frac{|h(x)-h(y)|}{\|x-y\|}$\vspace*{3pt}.
\item[(2)] For any $C^2$-function $h\dvtx \R^d\to\R$ we have
\[
\bigl|E \bigl[h(F) \bigr] - E \bigl[h(N) \bigr]\bigr | \leq\tfrac12\bigl\|h''
\bigr\|_\infty \sqrt{E\|F\|^4-E\|N\|^4},
\]
where $\|h''\|_\infty=
\max_{1\leq i,j\leq d}   \sup_{x\in\R^d}\llvert  \frac{\partial^2
h}{\partial x_i\,\partial x_j}(x)\rrvert  $.
\end{longlist}
\end{thmm}
\begin{pf}
The proof is divided into three steps.

\textit{Step} 1. Recall that for a Lipschitz function
$h\dvtx \R^d\to\R$, \cite{ivangiobook}, Theorem 6.1.1, yields\vspace*{1pt}
\begin{eqnarray*}
&&\bigl|E \bigl[h(F) \bigr] - E \bigl[h(N) \bigr] \bigr| \\
&&\qquad\leq\sqrt{d} \|\Sigma
\|_{\mathrm{op}}^{1/2}\bigl\|\Sigma^{-1}\bigr\|_{\mathrm{op}}\|h
\|_{\mathrm{Lip}}\sqrt {\sum_{i,j=1}^d
E \biggl\{ \biggl( \sigma_{ij}-\frac{1}{q_j}\langle
DF_i,DF_j\rangle \biggr)^2 \biggr\}},
\end{eqnarray*}
while for a $C^2$-function with bounded Hessian,
\cite{ivangiobook}, Theorem 6.1.2, gives\vspace*{1pt}
\begin{eqnarray*}
&&\bigl|E \bigl[h(F) \bigr] - E \bigl[h(N) \bigr] \bigr| \\
&&\qquad\leq\frac12\bigl\|h''
\bigr\|_\infty\sqrt{\sum_{i,j=1}^dE
\biggl\{ \biggl( \sigma_{ij}-\frac{1}{q_j}\langle
DF_i,DF_j\rangle \biggr)^2 \biggr\}}.
\end{eqnarray*}

\textit{Step} 2. We claim that for any $i,j=1,\ldots,d$,\vspace*{1pt}
\begin{eqnarray*}
&&E \biggl\{ \biggl( \sigma_{ij}-\frac{1}{q_j}\langle
DF_i,DF_j\rangle \biggr)^2 \biggr\}\\
&&\qquad\leq\operatorname{Cov }\bigl(F_i^2,F_j^2
\bigr) - 2\sigma_{ij}^2.
\end{eqnarray*}
Indeed, by \cite{ivangiobook}, identity (6.2.4), and the fact that
$\sigma_{ij}=0$ if $q_i\neq q_j$, we have\vspace*{1pt}
\begin{eqnarray*}
&&E \biggl\{ \biggl( \sigma_{ij}-\frac{1}{q_j}\langle
DF_i,DF_j\rangle \biggr)^2 \biggr\}
\\[3pt]
&&\qquad=\cases{ %
\displaystyle q_i^2\sum
_{r=1}^{q_i\wedge q_j}(r-1)!^2
\pmatrix{q_i-1
\cr
r-1}^2 \pmatrix{q_j-1
\cr
r-1}^2(q_i+q_j-2r)!\| f_i
\,\widetilde{\otimes}_r\, f_j\| ^2,\vspace*{4pt}\cr
\qquad\mbox{if
$q_i\neq q_j$},
\vspace*{4pt}\cr
\displaystyle q_i^2\sum_{r=1}^{q_i-1}(r-1)!^2
\pmatrix {q_i-1
\cr
r-1}^4 (2q_i-2r)!\|
f_i\,\widetilde{\otimes}_r\, f_j
\|^2,
\vspace*{4pt}\cr
\qquad\mbox{if $q_i= q_j$},}
\\[3pt]
&&\qquad\leq\cases{ %
\displaystyle\sum
_{r=1}^{q_i\wedge q_j}r!^2\pmatrix{q_i
\cr
r}^2 \pmatrix{q_j
\cr
r}^2(q_i+q_j-2r)!
\| f_i\,\widetilde{\otimes}_r\, f_j
\|^2,&\quad $\mbox {if $q_i\neq q_j$},$
\vspace*{4pt}\cr
\displaystyle\sum_{r=1}^{q_i-1}r!^2
\pmatrix{q_i
\cr
r}^4 (2q_i-2r)!\|
f_i\,\widetilde{\otimes}_r\, f_j
\|^2,&\quad $\mbox{if $q_i= q_j$}$}\
\end{eqnarray*}

On the other hand, from (\ref{cov}) we have
\begin{eqnarray*}
\hspace*{-4pt}&&\operatorname{Cov} \bigl(F_i^2,F_j^2
\bigr) - 2\sigma_{ij}^2
\\
\hspace*{-10pt}&&\qquad=\cases{ %
\displaystyle q_i!q_j!
\sum_{r=1}^{q_i\wedge q_j}\pmatrix{q_i
\cr
r} \pmatrix{q_j
\cr
r}\| f_i\otimes_r
f_j\|^2
\vspace*{2pt}\cr
\displaystyle\hspace*{-6pt}\qquad {}+\sum_{r=1}^{q_i\wedge q_j}r!^2
\pmatrix{q_i
\cr
r}^2 \pmatrix{q_j
\cr
r}^2(q_i+q_j-2r)!\| f_i
\,\widetilde{\otimes}_r\, f_j\|^2, \hspace*{17pt}\mbox{if
$q_i\neq q_j$},
\vspace*{2pt}\cr
\displaystyle q_i!^2\sum_{r=1}^{q_i-1}
\pmatrix{q_i
\cr
r}^2\| f_i
\otimes_r f_j\|^2
\vspace*{2pt}\cr
\displaystyle\hspace*{-6pt}\qquad{}+\sum_{r=1}^{q_i-1}r!^2
\pmatrix{q_i
\cr
r}^4 (2q_i-2r)!\|
f_i\,\widetilde{\otimes}_r\, f_j
\|^2,\hspace*{72pt}\mbox{if $q_i= q_j$},}
\end{eqnarray*}
The claim follows immediately.

\textit{Step} 3. Applying \eqref{identity} we get
\begin{eqnarray*}
E\|F\|^4-E\|N\|^4&=&\sum_{i,j=1}^d
\bigl(E \bigl[F_i^2F_j^2
\bigr] - \sigma_{ii} \sigma_{jj} - 2 \sigma_{ij}^2
\bigr)
\\
&=&\sum_{i,j=1}^d \bigl\{\operatorname{Cov}
\bigl(F_i^2,F_j^2 \bigr) -2
\sigma_{ij}^2 \bigr\}.
\end{eqnarray*}

Combining Steps 1--3 gives the desired conclusion.
\end{pf}

\subsection{A multivariate version of the convergence toward $\chi^2$}

Here we will prove a multivariate extension of a result of Nourdin and
Peccati~\cite{NPchi}. Such an extension was an open problem as far as
we know.

In what follows, $G(\nu)$ will denote a random variable with the
centered $\chi^2$ distribution having $\nu>0$ degrees of freedom.
When $\nu$ is an integer, then
$G(\nu)\stackrel{\mathrm{law}}{=}\sum_{i=1}^\nu(N_i^2-1)$, where $N_1,\ldots
,N_\nu$ are i.i.d. standard normal random variables. In general, $G(\nu
)$ is a centered gamma random variable with a shape parameter $\nu/2$
and scale parameter $2$.
Nourdin and Peccati~\cite{NPchi} established the following theorem.

\begin{thmm}[(Nourdin--Peccati)]\label{noupec}
Fix $\nu>0$, and let $G(\nu)$ be as above.
Let $q\geq2$ be an even integer, and
let $F_n=I_q(f_n)$ be such that $\lim_{n\to\infty} E[F_n^2] =E[G(\nu
)^2]=2\nu$.
Set $c_q= 4 [(q/2)! ]^3  [q! ]^{-2}$.
Then the following four assertions are equivalent, as $n\to\infty$:
\begin{longlist}[($\alpha$)]
\item[($\alpha$)] $F_n \stackrel{\mathrm{law}}{\to}G(\nu)$;
\item[($\beta$)] $E[F_n^4] -12E[F_n^3] \to E[G(\nu)^4] -12E[G(\nu)^3] =
12\nu^2 - 48\nu$;
\item[($\gamma$)] $\|f_n\,\widetilde{\otimes}_{q/2}\, f_n-c_q\times f_n\|
\to0$, and
$\|f_n \otimes_{r} f_n\|\to0$ for every $r=1,\ldots,q-1$ such that
$r\neq q/2$;
\item[($\delta$)] $\|f_n\,\widetilde{\otimes}_{q/2}\, f_n-c_q\times f_n\|
\to0$, and
$\|f_n\, \widetilde{\otimes}_r\, f_n\|\to0$ for every $r=1,\ldots,q-1$ such that
$r\neq q/2$.
\end{longlist}
\end{thmm}

The following is our multivariate extension of this theorem.

\begin{thmm}\label{tNouPec}
Let $d\geq2$, let $\nu_1,\ldots,\nu_d$ be positive reals and let
$q_1,\ldots,\break  q_d\geq2$ be
\textit{even} integers. Consider vectors
\[
F_n=(F_{1,n},\ldots,F_{d,n}) =
\bigl(I_{q_1}(f_{1,n}),\ldots, I_{q_d}(f_{d,n})
\bigr),\qquad n\geq1,
\]
with $f_{i,n}\in\HH^{\odot q_i}$, such that
$\lim_{n\to\infty} E[F_{i,n}^2] = 2\nu_i$ for every $i=1,\ldots,d$.
Assume that:
\begin{longlist}[(iii)]
\item[(i)] $E[F_{i,n}^4] -12E[F_{i,n}^3] \to12\nu_i^2 - 48\nu_i$  for
every $i$;

\item[(ii)]
$\lim_{n\to\infty} \operatorname{Cov}(F_{i,n}^2,F_{j,n}^2) = 0$
whenever $q_i=q_j$ for some $i\neq j$;

\item[(iii)]
$\lim_{n\to\infty} E[F_{i,n}^2F_{j,n}] = 0$ whenever $q_j=2q_i$.
\end{longlist}
Then
\[
(F_{1,n},\ldots,F_{d,n})\stackrel{\mathrm{law}} {\to} \bigl(G(
\nu_1),\ldots,G(\nu_d) \bigr),
\]
where $G(\nu_1),\ldots,G(\nu_d)$ are independent random variables
having centered $\chi^2$ distributions with $\nu_1,\ldots,\nu_d$
degrees of freedom, respectively.
\end{thmm}

%

\begin{pf}
Using the well-known Carleman condition, it is easy to check that the
law of $G(\nu)$ is determined by its moments.
By Corollary \ref{marginalblock} it is enough to show that condition
$(\gamma)$ of Theorem \ref{main} holds.

Fix $1\leq i\neq j\leq d$ as well as $1\leq r\leq q_i\wedge q_j$.
Switching $i$ and $j$ if necessary, assume that $q_i\leq q_j$.
From Theorem \ref{noupec}$(\gamma)$ we get that
$f_{k,n}\otimes_r f_{k,n}\to0$ for each $1\leq k\leq d$ and every
$1\leq r\leq q_k-1$, except when $r= q_k/2$.
Using the identity
%
\begin{equation}
\label{cs} \|f_{i,n}\otimes_r f_{j,n}
\|^2=\langle f_{i,n}\otimes_{q_i-r}
f_{i,n}, f_{j,n}\otimes_{q_j-r} f_{j,n}
\rangle
\end{equation}
[see (\ref{fubini})]
together with the Cauchy--Schwarz inequality, we infer that condition
$(\gamma)$ of Theorem \ref{main} holds for all values of $r$, $i$ and
$j$, except for the cases: $r=q_i=q_j$, $r=q_i/2=q_j/2$ and
$r=q_i=q_j/2$. Assumption (i) together with (\ref{cov2}) show that
$f_{i,n}\otimes_r f_{j,n}\to0$ for all
$1\leq r\leq q_i=q_j$. Thus it remains to verify condition $(\gamma)$
of Theorem \ref{main} when $r=q_i=q_j/2$.
Lemma \ref{technical} [identity (\ref{azertiop}) therein] yields
\begin{eqnarray*}
&&\langle f_{j,n}\, \widetilde{\otimes}_{q_i}\,
f_{j,n}, f_{i,n} \,\widetilde{\otimes}\, f_{i,n}\rangle
\\
&&\qquad= \frac{2q_i!^2}{q_j!} \langle f_{j,n}\otimes_{q_i}
f_{j,n},f_{i,n}\otimes f_{i,n}\rangle\\
&&\qquad\quad{} +
\frac{q_i!^2}{q_j!}\sum_{s=1}^{q_i-1}
\pmatrix{q_i
\cr
s}^2\langle f_{j,n}
\otimes_s f_{i,n},f_{i,n}\otimes_s
f_{j,n}\rangle.
\end{eqnarray*}
Using (\ref{cs}) and Theorem \ref{noupec} and reasoning as above, it is
straightforward to show that the sum $\sum_{s=1}^{q_i-1}{q_i\choose s}^2\langle f_{j,n}\otimes_s f_{i,n},f_{i,n}\otimes_s
f_{j,n}\rangle$
tends to zero as $n\to\infty$.
On the other hand,
the condition on the $q_i$th contraction in Theorem \ref{noupec}($\delta
$) yields that
$f_{j,n}\, \widetilde{\otimes}_{q_i}\,  f_{j,n} -
c_{q_j} f_{j,n}\to0$
as $n\to\infty$. Moreover, we have
\[
\langle f_{j,n}, f_{i,n} \,\widetilde{\otimes}\,
f_{i,n}\rangle =\frac{1}{q_j!}E \bigl[F_{j,n}F_{i,n}^2
\bigr],
\]
which tends to zero by assumption (ii).
All these facts together imply that
${\langle f_{j,n}\otimes_{q_i} f_{j,n},f_{i,n}\otimes f_{i,n}\rangle\to
0}$ as $n\to\infty$. Using (\ref{cs}) for $r=q_i$ we get $f_{i,n}
\otimes_{q_i} f_{j,n} \to0$, showing that condition $(\gamma)$ of
Theorem \ref{main} holds true in the last remaining case.
The proof of the theorem is complete.
\end{pf}

\begin{example}\label{}
Consider $F_n=(F_{1,n}, F_{2,n})=(I_{q_1}(f_{1,n}), I_{q_2}(f_{2,n}))$,
where $2\le q_1 \le q_2$ are even integers. Suppose that
\begin{eqnarray*}
E \bigl[F_{1,n}^2 \bigr] &\to& 1, \qquad E \bigl[F_{1,n}^4
\bigr] -6E \bigl[F_{1,n}^3 \bigr] \to-3 \quad\mbox{and}
\\
E \bigl[F_{2,n}^2 \bigr] &\to& 2,\qquad E \bigl[F_{2,n}^4
\bigr] -6E \bigl[F_{2,n}^3 \bigr] \to0\qquad \mbox{as } n\to
\infty.
\end{eqnarray*}
When $q_1=q_2$ or $q_2=2q_1$ we require additionally,
\[
\operatorname{Cov} \bigl(F_{1,n}^2,F_{2,n}^2
\bigr) \to0 \qquad(q_1=q_2), \qquad E \bigl[F_{1,n}^2F_{2,n}
\bigr] \to0 \qquad (q_2=2q_1).
\]
Then Theorem \ref{tNouPec} (the case $\nu_1=2$, $\nu_2=4$) gives
\[
F_n \stackrel{\mathrm{law}} {\to} (V_1-1,
V_2+V_3-2),
\]
where $V_1,V_2,V_3$ are i.i.d. standard exponential random variables.
\end{example}

\subsection{Bivariate convergence}

\begin{thmm}\label{tbivariateblock}
Let $p_1,\ldots,p_r,q_1,\ldots,q_s$ be positive integers.
Assume further that $\min p_i\geq\max q_j$.
Consider
\begin{eqnarray*}
&&(F_{1,n},\ldots,F_{r,n},G_{1,n},\ldots,
G_{s,n}) \\
&&\qquad= \bigl(I_{p_1}(f_{1,n}),
\ldots,I_{p_r}(f_{r,n}),I_{q_1}(g_{1,n}),
\ldots,I_{q_s}(g_{s,n}) \bigr),\qquad n\geq1,
\end{eqnarray*}
with $f_{i,n}\in\HH^{\odot p_i}$ and $g_{j,n}\in\HH^{\odot q_j}$.
Suppose that as $n\to\infty$
%
\begin{equation}
\label{marg} F_{n}=(F_{1,n},\ldots,F_{r,n})
\stackrel{\mathrm{law}} {\to} N \quad\mbox {and}\quad G_{n}=(G_{1,n},
\ldots, G_{s,n}) \stackrel{\mathrm{law}} {\to} V,
\end{equation}
where $N\sim N_r(0,\Sigma)$, the marginals of $V$ are determined by
their moments and $N, V$ are independent.
If $E[F_{i,n} G_{j,n}] \to0$ (which trivially holds when $p_i\ne q_j$)
for all $i,j$, then
%
\begin{equation}
\label{biv} (F_{n},G_{n})\stackrel{\mathrm{law}} {\to}(N,V)
\end{equation}
jointly, as $n \to\infty$.
\end{thmm}

\begin{pf}
We will show that condition $(\mathbf{c})$ of Theorem \ref
{mainblock} holds. By \eqref{bounded} we may and do assume that
$E[F_{i,n}^2]=1$ for all $i$ and $n$. By Theorem \ref{nupec}$(\gamma)$,
$\|f_{i,n}\otimes_r f_{i,n}\|\to0$ for all $r=1,\ldots,p_i-1$.
Observe that
\[
\|f_{i,n}\otimes_r g_{j,n}\|^2=
\langle f_{i,n}\otimes_{p_i-r} f_{i,n},
g_{j,n}\otimes_{q_j-r} g_{j,n}\rangle
\]
so that $\|f_{i,n}\otimes_r g_{j,n}\| \to0$ for $1\le r \le p_i\wedge
q_j =q_j$, except possibly when $r=p_i=q_j$. But in this latter case,
\[
p_i!\|f_{i,n}\otimes_r g_{j,n}\| =
p_i! \bigl|\langle f_{i,n},g_{j,n}\rangle \bigr| =
\bigl|E[F_{i,n}G_{j,n}] \bigr| \to0
\]
by the assumption. Corollary \ref{marginalblock} completes the proof.
\end{pf}

Theorem \ref{tbivariateblock} admits the following immediate corollary.

\begin{cor}\label{corindblocks}
Let $p\geq q$ be positive integers.
Consider two stochastic processes $F_n=(I_p(f_{t,n}))_{t\in T}$
and $G_n=(I_q(g_{t,n}))_{t\in T}$,
where $f_{t,n}\in\HH^{\odot p}$ and $g_{t,n}\in\HH^{\odot q}$. Suppose
that as $n\to\infty$,
\[
F_{n} \stackrel{{f.d.d.}} {\to} X \quad\mbox{and}\quad G_{n}
\stackrel{{f.d.d.}} {\to} Y,
\]
where $X$ is centered and Gaussian, the marginals of $Y$ are determined
by their moments and $X, Y$ are independent.
If $E[I_p(f_{t,n}) I_q(g_{s,n})] \to0$ (which trivially holds when
$p\ne q$) for all $s,t\in T$, then
\[
(F_{n},G_{n})\stackrel{{f.d.d.}} {\to}(X,Y)
\]
jointly, as $n \to\infty$.
\end{cor}

\section{Further applications}\label{sf-apps}

\subsection{Partial sums associated with Hermite polynomials}
Consider a centered stationary Gaussian sequence $\{G_k\}_{k\geq1}$
with unit variance.
For any $k\geq0$, denote by
\[
r(k)=E[G_1G_{1+k}]\vadjust{\goodbreak}
\]
the covariance between $G_1$ and $G_{1+k}$. We extend $r$ to $\Z_-$ by
symmetry, that is, $r(k)=r(-k)$.
For any integer $q\geq1$, we write
\[
S_{q,n}(t)= \sum_{k=1}^{\lfloor nt\rfloor}
H_q(G_k),\qquad t\geq0,
\]
to indicate the partial sums associated with the subordinated sequence
$\{H_q(G_k)\}_{k\geq1}$. Here, $H_q$ denotes the $q$th Hermite
polynomial given by (\ref{hermitepol}).

The following result is a summary of the main finding in Breuer and
Major~\cite{BM}.
%
\begin{thmm}\label{BM}
If $\sum_{k\in\Z} |r(k)|^q<\infty$, then as $n\to\infty$,
\[
\frac{S_{q,n}}{\sqrt{n}}\stackrel{{f.d.d.}} {\to} a_q B,
\]
where $B$ is a standard Brownian motion and $a_q= [q!\sum_{k\in\Z
}r(k)^q ]^{1/2}$.
\end{thmm}

Assume further that the covariance function $r$ has the form
\[
r(k)=k^{-D}L(k),\qquad k\geq1,
\]
with $D>0$ and $L\dvtx (0,\infty)\to(0,\infty)$ a function which is slowly
varying at infinity and bounded away from $0$ and infinity on every
compact subset of $[0,\infty)$.
The following result is due to Taqqu \cite{Taqqu}.

\begin{thmm}\label{taqqu}
If $0<D<\frac{1}{2}$, then as $n\to\infty$,
\[
\frac{S_{2,n}}{n^{1-D}L(n)} \stackrel{{f.d.d.}} {\to} b_D
R_{1-D},
\]
where $b_D= [(1-D)(1-2D) ]^{-1/2}$, and
$R_{H}$ is a Rosenblatt process of parameter $H=1-D$,
defined as
\[
R_{H}(t)=c_H I_2 \bigl(f_{H}(t,
\cdot) \bigr),\qquad t\geq0,
\]
with
\[
f_{H}(t,x,y)=\int_0^t
(s-x)_+^{{H}/{2}-1}(s-y)_+^{{H}/{2}-1}\,ds,\qquad t\geq0, x,y\in \R,
\]
$c_H>0$ an explicit constant such that $E[R_{H}(1)^2]=1$, and the
double Wiener--It\^o integral $I_2$ is with respect to a two-sided
Brownian motion $B$.
\end{thmm}

Let $q\geq3$ be an integer. The following result is a consequence of
Corollary \ref{corindblocks} and Theorems \ref{BM} and \ref{taqqu}. It
gives the asymptotic behavior (after proper renormalization of each
coordinate) of the pair $(S_{q,n},S_{2,n})$ when $D\in (\frac
{1}{q},\frac{1}{2} )\cup (\frac{1}{2},\infty)$.
Since what follows is just meant to be an illustration, we will not
consider the remaining case, that is, when
$D\in (0,\frac{1}{q} )$; it is an interesting problem, but to
answer it
would be out of the scope of the present paper.

\begin{prop}\label{BM-DMT}
Let $q\geq3$ be an integer, and let the constants $a_p$ and $b_D$ be
given by Theorems \ref{BM} and \ref{taqqu}, respectively.
\begin{longlist}[(1)]
\item[(1)] If $D\in(\frac12,\infty)$, then
\[
\biggl(\frac{S_{q,n}}{\sqrt{n}},\frac{S_{2,n}}{\sqrt{n}} \biggr) \stackrel{{f.d.d.}} {\to}
( a_q B_1, a_2 B_2 ),
\]
where $(B_1,B_2)$ is a standard Brownian motion in $\R^2$.
\item[(2)] If $D\in (\frac{1}{q},\frac{1}{2} )$, then
\[
\biggl(\frac{S_{q,n}}{\sqrt{n}},\frac{S_{2,n}}{n^{1-D}L(n)} \biggr)
\\
\stackrel{{f.d.d.}} {\to} (a_q B, b_D
R_{1-D} ),
\]
where $B$ is a Brownian motion independent of the Rosenblatt process
$R_{1-D}$ of parameter $1-D$.
\end{longlist}
\end{prop}

\begin{pf}
Let us first introduce a specific realization of the sequence $\{G_k\}
_{k\geq1}$ that will allow one to use the results of this paper. The space
\[
\mathcal{H}:=\overline{\operatorname{span}\{G_1,G_2,\ldots
\}}^{L^2(\Omega)}
\]
being a real separable Hilbert space, is isometrically isomorphic to
either $\R^N$ (for some finite $N\geq1$) or
$L^2(\R_+)$. Let us assume that $\mathcal{H}\simeq L^2(\R_+)$, the case
where $\mathcal{H}\simeq\R^N$ being easier to handle. Let $\Phi
\dvtx \mathcal{H}\to L^2(\R_+)$ be an isometry. Set $e_k=\Phi(G_k)$ for each
$k\geq1$. We have
%
\begin{equation}
r(k-l)=E[G_kG_l]=\int_0^\infty
e_k(x)e_l(x)\,dx,\qquad k,l\geq1.\label{ekrho}
\end{equation}
If $B=(B_t)_{t\in\R_+}$ denotes a standard Brownian motion, we deduce that
\[
\{G_k\}_{k\geq1} \stackrel{\mathrm{law}} {=} \biggl\{\int
_0^\infty e_k(t)\,dB_t
\biggr\}_{k\geq1},
\]
these two sequences being indeed centered, Gaussian and having the same
covariance structure. Using (\ref{mapping}) we deduce that $S_{q,n}$
has the same distribution as $I_q (\sum_{k=1}^n e_k^{\otimes q}
)$ (with $I_q$ the $q$-tuple Wiener--It\^o integral associated to $B$).

Hence, to reach the conclusion of point 1 it suffices to combine
Corollary \ref{corindblocks} with Theorem \ref{BM}. For point 2,
just use
Corollary \ref{corindblocks} and Theorem \ref{taqqu},
together with the fact that the distribution of $R_H(t)$ is determined
by its moments (as is the case for any double Wiener--It\^o integral).
\end{pf}

\subsection{Moment-independence for discrete homogeneous chaos}
To develop the next application we will need the following basic ingredients:
\begin{enumerate}[(ii)]
\item[(i)] A sequence $\mathbf{X}=(X_1,X_2,\ldots)$ of i.i.d. random variables,
with mean 0, variance 1 and all moments finite.\vadjust{\goodbreak}
\item[(ii)] Two positive integers $q_1,q_2$ as well as two sequences
$a_{k,n}\dvtx \N^{q_k}\to\R$, $n\geq1$ of real-valued functions satisfying
for all $i_1,\ldots,i_{q_k} \ge1$ and $k=1,2$:
\begin{enumerate}[(a)]
\item[(a)] (symmetry) $a_{k,n}(i_1,\ldots,i_{q_k})=a_{k,n}(i_{\sigma
(1)},\ldots,i_{\sigma(q_k)})$ for every permutation $\sigma$;
\item[(b)] (vanishing on diagonals) $a_{k,n}(i_1,\ldots,i_{q_k})=0$
whenever $i_r=i_s$ for some $r\neq s$;
\item[(c)] (unit-variance)  $q_k!\sum_{i_1,\ldots,i_{q_k}=1}^{\infty}
a_{k,n}(i_1,\ldots,i_{q_k})^2=1$.
\end{enumerate}
\end{enumerate}
Consider
%
\begin{equation}\label{}
\qquad Q_{k,n}(\mathbf{X})=\sum_{i_1,\ldots,i_{q_k}=1}^{\infty}
a_{k,n}(i_1,\ldots,i_{q_k})X_{i_1}
\cdots X_{i_{q_k}}, \qquad n\geq1, k=1,2.
\end{equation}
This series converges in $L^2(\Omega)$, $E [Q_{k,n}(\mathbf{X})]=0$ and $E
[Q_{k,n}(\mathbf{X})^2]=1$.
We have the following result.


\begin{thmm}\label{thmMOO}
As $n\to\infty$, assume that the contribution of each $X_i$ to
$Q_{k,n}(\mathbf{X})$ is uniformly negligible, that is,
%
\begin{equation}
\sup_{i \ge1} \sum_{i_2,\ldots,i_{q_k}=1}^{\infty}
a_{k,n}(i,i_2,\ldots,i_{q_k})^2\to0,\qquad
k=1,2, \label{influence}
\end{equation}
and that, for any $r=1,\ldots,q_1\wedge q_2$,
%
\begin{eqnarray}\label{mixed}\qquad
&&\sum_{i_1,\ldots,i_{q_1+q_2-2r}=1}^{\infty} \Biggl(\sum
_{l_1,\ldots
,l_r=1}^{\infty}
 a_{1,n}(l_1,
\ldots,l_r,i_1,\ldots,i_{q_1-r})\nonumber
\\[-8pt]
\\[-8pt]
\nonumber
&&\hspace*{105pt}{}\times a_{2,n}(l_1,
\ldots,l_{r},i_{q_1-r+1},\ldots,i_{q_1+q_2-2r})
\Biggr)^2\to0.
\end{eqnarray}
Then
$Q_{1,n}(\mathbf{X})$ and $Q_{2,n}(\mathbf{X})$ are asymptotically moment-independent.
\end{thmm}

\begin{pf}
Fix $M,N\geq1$. We want to prove that, as
$n\to\infty$,
%
\begin{equation}
\label{goal} E \bigl[Q_{1,n}(\mathbf{X})^MQ_{2,n}(
\mathbf{X})^N \bigr]-E \bigl[Q_{1,n}(\mathbf{
X})^M \bigr]E \bigl[Q_{2,n}(\mathbf{X})^N
\bigr] \to0.
\end{equation}
The proof is divided into three steps.

\textit{Step} 1. In this step we show that
%
\begin{equation}
\label{MOOcsq} \qquad E \bigl[Q_{1,n}(\mathbf{X})^MQ_{2,n}(
\mathbf{X})^N \bigr]-E \bigl[Q_{1,n}(\mathbf{
G})^MQ_{2,n}(\mathbf{G})^N \bigr]\to0
\qquad\mbox{as } n \to\infty.
\end{equation}
Following the approach of Mossel, O'Donnel and Oleszkiewicz \cite{MOO},
we will use the Lindeberg replacement trick. Let
$\mathbf{G}=(G_1,G_2,\ldots)$ be a sequence of i.i.d. $N(0,1)$ random
variables independent of $\mathbf{X}$.
For a positive integer $s$, set $\mathbf{W}^{(s)}=(G_1,\ldots,G_s,X_{s+1},
X_{s+2},\ldots)$, and put $\mathbf{W}^{(0)}=\mathbf{X}$.
Fix $s \ge1$ and write for $k=1,2$ and $n \ge1$,
\begin{eqnarray*}
U_{k,n,s}&=&\sum_{ \stackrel{i_1,\ldots,i_{q_k}}{i_1\neq s, \ldots,
i_{q_k}\neq s}
} a_{k,n}(i_1,
\ldots,i_{q_k})W_{i_1}^{(s)}\cdots
W_{i_{q_k}}^{(s)},
\\
V_{k,n,s}&=&\sum_{ \stackrel{ i_1,\ldots,i_{q_k}}{\exists j\dvtx i_j= s} } a_{k,n}(i_1,
\ldots,i_{q_k})W_{i_1}^{(s)}\cdots\widehat{W}_{s}^{(s)}\cdots W_{i_{q_k}}^{(s)}
\\
&=&q_k\sum_{i_2,\ldots,i_{q_k}=1}^{\infty}
a_{k,n}(s, i_2,\ldots,i_{q_k})W_{i_2}^{(s)}
\cdots W_{i_{q_k}}^{(s)},
\end{eqnarray*}
where $\widehat{W}_s^{(s)}$ means that the term $W_s^{(s)}$ is dropped
(observe that this notation bears no ambiguity: indeed, since $a_{k,n}$
vanishes on diagonals,
each string $i_1,\ldots,i_{q_k}$ contributing to the definition of
$V_{k,n,s}$ contains the symbol
$s$ exactly once).
For each $s$ and $k$, note that $U_{k,n,s}$ and $V_{k,n,s}$ are
independent of
the variables $X_s$ and $G_s$, and that
\[
Q_{k,n} \bigl(\mathbf{W}^{(s-1)} \bigr)=U_{k,n,s}+X_sV_{k,n,s}
\quad\mbox{and}\quad Q_{k,n} \bigl(\mathbf{W}^{(s)}
\bigr)=U_{k,n,s}+G_sV_{k,n,s}.
\]
By the binomial formula, using the independence of $X_s$ from
$U_{k,n,s}$ and $V_{k,n,s}$, we have
\begin{eqnarray*}
&& E \bigl[Q_{1,n} \bigl(\mathbf{W}^{(s-1)} \bigr)^M
Q_{2,n} \bigl(\mathbf{W}^{(s-1)} \bigr)^N \bigr]
\\
&&\qquad=\sum_{i=0}^M\sum
_{j=0}^N \pmatrix{M
\cr
i}\pmatrix{N
\cr
j}E
\bigl[U_{1,n,s}^{M-i}U_{2,n,s}^{N-j}V_{1,n,s}^{i}V_{2,n,s}^j
\bigr]E \bigl[X_s^{i+j} \bigr].
\end{eqnarray*}
Similarly,
\begin{eqnarray*}
&&E \bigl[Q_{1,n} \bigl(\mathbf{W}^{(s)} \bigr)^M
Q_{2,n} \bigl(\mathbf{W}^{(s)} \bigr)^N \bigr]
\\
&&\qquad=\sum_{i=0}^M\sum
_{j=0}^N \pmatrix{M
\cr
i}\pmatrix{N
\cr
j}E
\bigl[U_{1,n,s}^{M-i}U_{2,n,s}^{N-j}V_{1,n,s}^{i}V_{2,n,s}^j
\bigr]E \bigl[G_s^{i+j} \bigr].
\end{eqnarray*}
Therefore
\begin{eqnarray*}
&&E \bigl[Q_{1,n} \bigl(\mathbf{W}^{(s-1)} \bigr)^M
Q_{2,n} \bigl(\mathbf{W}^{(s-1)} \bigr)^N \bigr] -
E \bigl[Q_{1,n} \bigl(\mathbf{W}^{(s)} \bigr)^M
Q_{2,n} \bigl(\mathbf{W}^{(s)} \bigr)^N \bigr]
\\
&&\qquad=\sum_{i+j \ge3} \pmatrix{M
\cr
i}\pmatrix{N
\cr
j} E
\bigl[U_{1,n,s}^{M-i}U_{2,n,s}^{N-j}V_{1,n,s}^{i}V_{2,n,s}^j
\bigr] \bigl( E \bigl[X_s^{i+j} \bigr] - E
\bigl[G_s^{i+j} \bigr] \bigr).
\end{eqnarray*}
Now, observe that Propositions 3.11, 3.12 and 3.16 of \cite{MOO} imply
that both
$(U_{1,n,s})_{n,s\ge1}$ and $(U_{2,n,s})_{n,s\ge1}$ are uniformly
bounded in all $L^p(\Omega)$ spaces.
It also implies that, for any $p\geq3$, $k=1,2$ and $n, s \ge1$,
\[
E \bigl[|V_{k,n,s}|^p \bigr]^{1/p} \leq
C_p E \bigl[V_{k,n,s}^2 \bigr]^{1/2},
\]
where $C_p$ depends only on $p$. Hence, for $0\le i \le M$, $0\le j \le
N$, $i+j\ge3$, we have
%
\begin{equation}
\label{d-bound} \bigl| E \bigl[U_{1,n,s}^{M-i}U_{2,n,s}^{N-j}
V_{1,n,s}^{i}V_{2,n,s}^j
\bigr] \bigr| \le C E \bigl[V_{1,n,s}^2 \bigr]^{i/2} E
\bigl[V_{2,n,s}^2 \bigr]^{j/2},
\end{equation}
where $C$ does not depend on $n,s\ge1$. Since $E[X_i]=E[G_i]=0$ and
$E[X_i^2]=E[G_i^2]=1$, we get
\[
E \bigl[V_{k,n,s}^2 \bigr] = q_kq_k!
\sum_{i_2,\ldots,i_{q_k}=1}^{\infty} a_{k,n}(s,i_2,
\ldots,i_{q_k})^2.
\]
When $i\ge3$, then \eqref{d-bound} is bounded from above by
\[
C \Biggl(\sup_{i\ge1} \sum_{i_2,\ldots,i_{q_1}=1}^{\infty}
a_{1,n}(i,i_2,\ldots,i_{q_1})^2
\Biggr)^{(i-2)/2} \sum_{i_2,\ldots
,i_{q_1}=1}^{\infty}
a_{1,n}(s,i_2,\ldots,i_{q_1})^2,
\]
where $C$ does not depend on $n,s\ge1$, and we get a similar bound
when $j\ge3$. If $i=2$, then $j\ge1$ ($i+j\ge3$), so \eqref{d-bound}
is bounded from above by
\[
C \Biggl(\sup_{i\ge1} \sum_{i_2,\ldots,i_{q_2}=1}^{\infty}
a_{2,n}(i,i_2,\ldots,i_{q_2})^2
\Biggr)^{j/2} \sum_{i_2,\ldots
,i_{q_1}=1}^{\infty}
a_{1,n}(s,i_2,\ldots,i_{q_1})^2,
\]
and we have a similar bound when $j=2$. Taking into account assumption
\eqref{influence} we infer that the upper-bound for \eqref{d-bound} is
of the form
\[
C \varepsilon_n \sum_{k=1}^2
\sum_{i_2,\ldots,i_{q_k}=1}^{\infty} a_{k,n}(s,i_2,
\ldots,i_{q_k})^2,
\]
where $\lim_{n \to\infty} \varepsilon_n=0$ and $C$ is independent of
$n,s$. We conclude that
\begin{eqnarray*}
&& \bigl| E \bigl[Q_{1,n} \bigl(\mathbf{W}^{(s-1)}
\bigr)^M Q_{2,n} \bigl(\mathbf{W}^{(s-1)}
\bigr)^N \bigr] - E \bigl[Q_{1,n} \bigl(
\mathbf{W}^{(s)} \bigr)^M Q_{2,n} \bigl(
\mathbf{W}^{(s)} \bigr)^N \bigr]\bigr |
\\
&&\quad \le C \varepsilon_n \sum_{k=1}^2
\sum_{i_2,\ldots
,i_{q_k}=1}^{\infty} a_{k,n}(s,i_2,
\ldots,i_{q_k})^2,
\end{eqnarray*}
where $C$ does not depend on $n,s$. Since, for fixed $k,n$,
$Q_{k,n}(\mathbf{W}^{(s)}) \to Q_{k,n}(\mathbf{G}) $ in $L^2(\Omega)$ as $s
\to\infty$, by Propositions 3.11, 3.12 and 3.16 of \cite{MOO}, the
convergence holds in all $L^p(\Omega)$. Hence
\begin{eqnarray*}
&&\bigl| E \bigl[Q_{1,n}(\mathbf{X})^M  Q_{2,n}(
\mathbf{X})^N \bigr]-E \bigl[Q_{1,n}(\mathbf{
G})^MQ_{2,n}(\mathbf{G})^N \bigr] \bigr|
\\
&&\quad\le\sum_{s=1}^{\infty} | E
\bigl[Q_{1,n} \bigl(\mathbf{W}^{(s-1)} \bigr)^M
Q_{2,n} \bigl(\mathbf{W}^{(s-1)} \bigr)^N \bigr] -
E \bigl[Q_{1,n} \bigl(\mathbf{W}^{(s)} \bigr)^M
Q_{2,n} \bigl(\mathbf{ W}^{(s)} \bigr)^N \bigr] |
\\
&&\quad \le C \varepsilon_n \sum_{k=1}^2
\sum_{i_1,\ldots,i_{q_k}=1}^{\infty} a_{k,n}(i_1,i_2,
\ldots,i_{q_k})^2 = C \bigl((q_1!)^{-1}+
(q_2!)^{-1} \bigr) \varepsilon_n.
\end{eqnarray*}
This proves \eqref{MOOcsq}.\vadjust{\goodbreak}

\textit{Step} 2.
We show that $n \to\infty$,
%
\begin{eqnarray}
\label{MOOcsq2} E \bigl[Q_{1,n}(\mathbf{X})^M \bigr]-E
\bigl[Q_{1,n}(\mathbf{G})^M \bigr]&\to&0 \quad\mbox{and}
\nonumber
\\[-8pt]
\\[-8pt]
\nonumber
 E
\bigl[Q_{2,n}(\mathbf{X})^N \bigr]-E \bigl[Q_{2,n}(
\mathbf{G})^N \bigr]&\to&0.
\end{eqnarray}
The proof is similar to Step 1 (and easier). Thus, we omit it.

\textit{Step} 3. Without loss of generality we may and do assume that
$G_k=B_k-B_{k-1}$, where $B$ is a standard Brownian motion. For $k=1,2$
and $n\geq1$, due to the multiplication formula (\ref
{multiplication}), $Q_{k,n}(\mathbf{G})$ is a multiple Wiener--It\^o
integral of order $q_k$ with respect to $B$,
\[
Q_{k,n}(\mathbf{G})=I_{q_k} \Biggl( \sum
_{i_1,\ldots,i_{q_k}=1}^{\infty} a_{k,n}(i_1,
\ldots,i_{q_k})\mathbf{ 1}_{[i_1-1,i_1]
\times\cdots\times[i_{q_k}-1,i_{q_k}]} \Biggr).
\]
In this setting, condition (\ref{mixed}) coincides with condition
($\gamma$) of Theorem \ref{main} [or ($\mathbf{c}$) of Theorem \ref
{mainblock}]. Therefore,
%
\begin{equation}
\label{mixedcsq} E \bigl[Q_{1,n}(\mathbf{G})^MQ_{2,n}(
\mathbf{G})^N \bigr]-E \bigl[Q_{1,n}(\mathbf{
G})^M \bigr]E \bigl[Q_{2,n}(\mathbf{G})^N
\bigr] \to0.
\end{equation}
Combining (\ref{MOOcsq}), (\ref{MOOcsq2}) and (\ref{mixedcsq}) we get
the desired conclusion (\ref{goal}).
\end{pf}

\begin{remark}
The conclusion of Theorem \ref{thmMOO} may fail if either (\ref
{influence}) or (\ref{mixed}) are not satisfied. It follows from Step 3
above that the theorem fails when \eqref{mixed} does not hold and $\mathbf{
X}$ is Gaussian. Theorem \ref{thmMOO} also fails when (\ref{influence})
is not satisfied, (\ref{mixed}) holds and $\mathbf{X}$ is a Rademacher
sequence, as we can see from the following counterexample. Consider
$q_1=q_2=2$, and set
\begin{eqnarray*}
a_{1,n}(i,j)&=&\tfrac14 \bigl(\mathbf{1}_{\{1\}}(i)
\mathbf{1}_{\{2\}}(j)+\mathbf{ 1}_{\{2\}}(i)\mathbf{1}_{\{1\}}(j)+
\mathbf{1}_{\{1\}}(i)\mathbf{1}_{\{3\}
}(j)+\mathbf{1}_{\{3\}}(i)
\mathbf{1}_{\{1\}}(j) \bigr),
\\
a_{2,n}(i,j)&=&\tfrac14 \bigl(\mathbf{1}_{\{2\}}(i)
\mathbf{1}_{\{4\}}(j)+\mathbf{ 1}_{\{4\}}(i)\mathbf{1}_{\{2\}}(j)-
\mathbf{1}_{\{3\}}(i)\mathbf{1}_{\{4\}
}(j)-\mathbf{1}_{\{4\}}(i)
\mathbf{1}_{\{3\}}(j) \bigr).
\end{eqnarray*}
Then
$Q_{1,n}(\mathbf{X})=\frac12 X_1(X_2+X_3)$ and $Q_{2,n}(\mathbf{X})=\frac12
X_4(X_2-X_3)$, where $X_i$ are i.i.d. with $P(X_i=1)=P(X_i=-1)=1/2$.
It is straightforward to check that~(\ref{mixed}) holds and obviously
(\ref{influence}) is not satisfied.
Since $Q_{1,n}(\mathbf{X})Q_{2,n}(\mathbf{X})=0$, we get
\[
0=E \bigl[Q_{1,n}(\mathbf{X})^2Q_{2,n}(
\mathbf{X})^2 \bigr]\neq E \bigl[Q_{1,n}(\mathbf{
X})^2 \bigr]E \bigl[Q_{2,n}(\mathbf{X})^2
\bigr],
\]
implying in particular that $Q_{1,n}(\mathbf{X})$ and $Q_{2,n}(\mathbf{X})$
are (asymptotically) moment-dependent.
\end{remark}

\section*{Acknowledgments}
We are grateful to Jean Bertoin for useful discussions and to
Ren\'e Schilling for reference \cite{BS}. We warmly thank an anonymous
referee for suggesting a shorter proof of Theorem \ref{main} (which
evolved into the proof of a more general statement, Theorem \ref{mainblock})
and for useful comments and suggestions, which together with the
Editor's constructive remarks, have led to a significant improvement of
this paper.




\printaddresses

\end{document}